\documentclass[reqno]{amsart}

\usepackage{epsfig}
\usepackage[all]{xy}

\def\E{\ifmmode{\mathbb E}\else{$\mathbb E$}\fi} 
\def\N{\ifmmode{\mathbb N}\else{$\mathbb N$}\fi} 
\def\R{\ifmmode{\mathbb R}\else{$\mathbb R$}\fi} 
\def\Q{\ifmmode{\mathbb Q}\else{$\mathbb Q$}\fi} 
\def\C{\ifmmode{\mathbb C}\else{$\mathbb C$}\fi} 
\def\H{\ifmmode{\mathbb H}\else{$\mathbb H$}\fi} 
\def\Z{\ifmmode{\mathbb Z}\else{$\mathbb Z$}\fi} 
\def\P{\ifmmode{\mathbb P}\else{$\mathbb P$}\fi} 
\def\T{\ifmmode{\mathbb T}\else{$\mathbb T$}\fi} 
\def\SS{\ifmmode{\mathbb S}\else{$\mathbb S$}\fi} 
\def\DD{\ifmmode{\mathbb D}\else{$\mathbb D$}\fi} 

\newcommand{\e}{\varepsilon}

\newcommand{\del}{\partial}

\newcommand{\ben}{\begin{enumerate}}
\newcommand{\een}{\end{enumerate}}
\newcommand{\be}{\begin{equation}}
\newcommand{\ee}{\end{equation}}
\newcommand{\bea}{\begin{eqnarray}}
\newcommand{\eea}{\end{eqnarray}}
\newcommand{\beastar}{\begin{eqnarray*}}
\newcommand{\eeastar}{\end{eqnarray*}}
\newcommand{\bc}{\begin{center}}
\newcommand{\ec}{\end{center}}

\theoremstyle{theorem}
\newtheorem{thm}{Theorem}[section]
\newtheorem{cor}[thm]{Corollary}
\newtheorem{lem}[thm]{Lemma}
\newtheorem{prop}[thm]{Proposition}

\theoremstyle{definition}
\newtheorem{defn}[thm]{Definition}

\newtheorem*{thm*}{Theorem}

\numberwithin{equation}{section}

\def\R{{\mathbb R}}
\def\osc{{\hbox{\rm osc}}}

\def\E{{\mathbb E}}
\def\Z{{\mathbb Z}}
\def\C{{\mathbb C}}
\def\R{{\mathbb R}}
\def\P{{\mathbb P}}

\def\N{{\mathbb N}}

\def\11{{\mathbb I}}

\def\C{\mathbb{C}}
\def\Z{\mathbb{Z}}

\def\T{\mathbb{T}}

\def\Q{\mathbb{Q}}

\def\E{\ifmmode{\mathbb E}\else{$\mathbb E$}\fi} 
\def\N{\ifmmode{\mathbb N}\else{$\mathbb N$}\fi} 
\def\R{\ifmmode{\mathbb R}\else{$\mathbb R$}\fi} 
\def\Q{\ifmmode{\mathbb Q}\else{$\mathbb Q$}\fi} 
\def\C{\ifmmode{\mathbb C}\else{$\mathbb C$}\fi} 
\def\H{\ifmmode{\mathbb H}\else{$\mathbb H$}\fi} 
\def\Z{\ifmmode{\mathbb Z}\else{$\mathbb Z$}\fi} 
\def\P{\ifmmode{\mathbb P}\else{$\mathbb P$}\fi} 
\def\SS{\ifmmode{\mathbb S}\else{$\mathbb S$}\fi} 
\def\DD{\ifmmode{\mathbb D}\else{$\mathbb D$}\fi} 

\def\R{{\mathbb R}}
\def\osc{{\hbox{\rm osc}}}

\def\E{{\mathbb E}}
\def\Z{{\mathbb Z}}
\def\C{{\mathbb C}}
\def\R{{\mathbb R}}

\def\N{{\mathbb N}}

\def\e{\varepsilon}

\def\CH{{\mathcal H}}

\def\CP{{\mathcal P}}
\def\CQ{{\mathcal Q}}

\def\CP{{\mathcal P}}

\def\darr#1{\raise1.5ex\hbox{$\leftrightarrow$}
\mkern-16.5mu #1}

\def\roughly#1{\raise.3ex\hbox{$#1$\kern-.75em
\lower1ex\hbox{$\sim$}}}

\def\opname#1{\mathop{\kern0pt{\rm #1}}\nolimits}

\begin{document}
\quad \vskip1.375truein

\def\mq{\mathfrak{q}}
\def\mp{\mathfrak{p}}
\def\mH{\mathfrak{H}}
\def\mh{\mathfrak{h}}
\def\ma{\mathfrak{a}}
\def\ms{\mathfrak{s}}
\def\mm{\mathfrak{m}}
\def\mn{\mathfrak{n}}

\def\Hoch{{\tt Hoch}}
\def\mt{\mathfrak{t}}
\def\ml{\mathfrak{l}}
\def\mT{\mathfrak{T}}
\def\mL{\mathfrak{L}}
\def\mg{\mathfrak{g}}
\def\md{\mathfrak{d}}

\newcommand{\includefigure}[3]{%
\begin{figure}[#3]
\begin{center}
\epsfig{file=#2} \\ \caption{\label{fig:#1}}
\end{center}
\end{figure}}

\title[Locality of continuous Hamiltonian flows]{Locality of
continuous Hamiltonian flows \\and Lagrangian intersections with\\
the conormal of open subsets}

\author{Yong-Geun Oh}
\thanks{Partially supported by the NSF grant \# DMS 0503954}
\date{Revision, January 2007}

\address{
Department of Mathematics, University of Wisconsin, Madison, WI
53706 \& Korea Institute for Advanced Study, 207-43
Cheongryangni-dong Dongdaemun-gu, Seoul 130-012, KOREA,
oh@math.wisc.edu}

\begin{abstract} In this paper, we prove that if a continuous
Hamiltonian flow fixes the points in an open subset $U$ of a
symplectic manifold $(M,\omega)$, then its associated Hamiltonian is
constant at each moment on $U$. As a corollary, we prove that the
Hamiltonian of compactly supported continuous Hamiltonian flows is
unique both on a compact $M$ with contact-type boundary $\del M$ and on a
non-compact manifold bounded at infinity. An essential tool for the
proof of the locality is the Lagrangian intersection theorem for the
conormals of open subsets proven by Kasturirangan and the author,
combined with Viterbo's scheme that he introduced in the proof of
uniqueness of the Hamiltonian on a closed manifold \cite{viterbo2}.
We also prove the converse of the theorem which localizes a
previously known global result in symplectic topology.
\end{abstract}


\maketitle

\tableofcontents

\section{Introduction}
\label{sec:intro}

Let $\phi$ be a Hamiltonian diffeomorphism on a symplectic
manifold $(M,\omega)$. Hofer's $L^\infty$-norm of Hamiltonian
diffeomorphisms is defined by
$$
\|\phi\| = \inf_{H\mapsto \phi} \|H\|_\infty
$$
where $H \mapsto \phi$ means that $\phi= \phi_H^1$ is the time-one
map of Hamilton's equation
$$
\dot x = X_H(t,x)
$$
and the norm $\|H\|$ is  defined by
\be\label{eq:intosc}
\|H\|_\infty = \max_{t \in [0,1]}\mbox{\rm osc }H_t
\ee
where
$$
\mbox{\rm osc }H_t = \max_x H_t -\min_x H_t.
$$
For two given
continuous paths $\lambda, \, \mu:[a,b] \to Homeo(M)$, we define
their distance by \be\label{eq:bard} \overline d(\lambda,\mu) =
\max_{t\in [a,b]}\overline d(\lambda(t),\mu(t)). \ee
\begin{defn}[Continuous Hamiltonian flow]\label{topflowdefn} A continuous map
$\lambda: \R \to Homeo(M)$ is called a \emph{continuous Hamiltonian flow}
if there exists a sequence of smooth Hamiltonians $H_i : \R \times
M \to \R$ satisfying the following :
\begin{enumerate}
\item $\phi_{H_i} \to \lambda$ locally uniformly on $\R \times M$.
\item the sequence $H_i$ is Cauchy in the $C^0$-topology and so
has a limit $H$ lying in $C^0$.
\end{enumerate}
We call a continuous path $\lambda:[a,b] \to Homeo(M)$ a {\it
continuous Hamiltonian path} if it satisfies the same conditions
with $\R$ replaced by $[a,b]$, and the limit $C^0$-function
$H$ a \emph{continuous Hamiltonian}. In any of these cases,
we say that the pair $(\lambda,H)$ is the \emph{$C^0$-hamiltonian
limit} of $(\phi_{H_i},H_i)$, and write
$$
\operatorname{hlim}_{i \to \infty} (\phi_{H_i},H_i) \to (\lambda,
H)
$$
or sometimes even $\operatorname{hlim}_{i \to \infty}
(\phi_{H_i},H_i)=\lambda$.
\end{defn}

We denote by $\CP^{ham}_{\infty;[a,b]}(Sympeo(M,\omega),id)$ the set of
continuous Hamiltonian paths defined on $[a,b]$. When $[a,b]=[0,1]$
or when we do not specify the domain of $\lambda$, we often just
write $\CP^{ham}_\infty(Sympeo(M,\omega),id)$ for the corresponding
set of continuous Hamiltonian paths.

For the rest of the paper until the last section, we will assume
that $M$ is closed.

An proof of the following theorem was provided by Viterbo in
\cite{viterbo2} for the $C^0$-hamiltonian limits. An immediate
corollary of the theorem is that each continuous Hamiltonian path
carries a unique Hamiltonian. This uniqueness question was asked by
the author in \cite{oh:hameo1}.

\begin{thm}[Viterbo, \cite{viterbo2}]
\label{uniqueness} Let $M$ be a closed symplectic manifold.
Suppose $H_i:[0,1] \times M \to \R$ is a sequence of smooth
normalized Hamiltonian functions converging to a
continuous function $H: Y \to \R$ in $C^0$-topology.
Then if $\phi_{H_i} \to id$ uniformly on $[0,1] \times M$, we have
$H \equiv 0$.
\end{thm}

The main purpose of the present paper is to prove the following
local version of Theorem \ref{uniqueness} and its converse.

\begin{thm}[$C^0$-Locality]\label{local-uniqueness}
Let $U \subset M$ be an open subset. Suppose $H_i:[0,1] \times M \to
\R$ is a sequence of smooth normalized Hamiltonian functions
such that both $\phi_{H_i}$ and $H_i$ converge in $C^0$-topology.
Denote their limits by $\lambda$ and $H$ respectively.
Then the followings are equivalent
\begin{enumerate}
\item $\lambda_t(x) = x$ for all $(t,x) \in [0,1] \times U$,
\item $H(t,x) \equiv c(t)$ on $U$ for a continuous function
$c: [0,1] \to \R$ depending only on $t$.
\end{enumerate}
\end{thm}

This locality of Hamiltonians is an easy and well-known fact for the
case of smooth (or more generally $C^1$) Hamiltonian paths. The case
$U = M$ is precisely the uniqueness question asked in
\cite{oh:hameo1} and answered by Viterbo. Our proof uses the crucial
idea from \cite{viterbo2} which reduces the uniqueness proof to a
certain Lagrangian intersection theorem.

One rather immediate consequence of the locality theorem above is the
following uniqueness theorem for \emph{compactly supported}
continuous Hamiltonian flows. We refer to \cite{oh:hameo1} for the
precise definition of compactly supported topological Hamiltonian
flows (i.e., $L^{(1,\infty)}$-case) which can be adapted to the
current continuous Hamiltonian flows (i.e., $C^0$-case).

\begin{thm}[$C^0$-Uniqueness on open manifolds]\label{compactsupp}
Let $M$ be a compact symplectic manifold with contact-type boundary $\del
M$ or a noncompact manifold bounded at infinity. Suppose $H_i:[0,1]
\times M \to \R$ is a sequence of smooth normalized Hamiltonian
functions such that
\begin{enumerate}
\item there exists a compact set $K \subset \operatorname{Int}M$
such that
$$
\operatorname{supp} H_i \subset K
$$
\item $H_i$ uniformly converges to a continuous function $H: Y \to \R$.
\end{enumerate}
Then if $\phi_{H_i} \to id$ uniformly on $[0,1] \times M$, we have
$H \equiv 0$.
\end{thm}

The special case of two-dimensional compact surface with smooth
boundary is an immediate consequence of the uniqueness result on
closed surfaces by the (space)-doubling argument (See
\cite{oh:calabi} for such a proof.) But this doubling argument
cannot be applied to the high dimensional cases. The presence of the
geometric conditions `with contact-type boundary' or
`bounded an infinity' in the hypothesis is an
artifact of our proof that relies on the Floer homology theory of
conormals from \cite{kasturi-oh}. It is not obvious to the author
whether the same uniqueness holds on arbitrary open symplectic
manifolds which are not necessarily bounded at infinity.

The main idea of our proof of the locality theorem is again to
reduce its proof to a version of Lagrangian intersection theorem.
Viterbo used the global version of Lagrangian intersection theorem
on the cotangent bundle from \cite{hofer1}, \cite{laud-sik},
\cite{gromov:pseudo}. For the current locality proof, the Lagrangian
intersection theorem that we use is the \emph{relative} version of
Lagrangian intersection theorem that Kasturirangan and the author
proved in \cite{kasturi-oh} using the study of Floer homology of
conormals of open subsets :

\begin{thm}[Kasturirangan-Oh, \cite{kasturi-oh}]\label{kasturi-oh}
Let $U\subset N$ be any open subset with smooth
boundary $\partial U$, and $H:T^*N \times [0,1] \to \R$ be any
smooth Hamiltonian. Then
$$
\phi^1_H(o_N)\cap \nu^*\overline{U} \neq \emptyset
$$
for all $t \in [0,1]$
where $\nu^*\overline U$ is the conormal of $\overline U$.
\end{thm}

We refer readers to section \ref{subsec:conormal} for the definition
of the conormal $\nu^*\overline U$. Once we have this local
Lagrangian intersection theorem, the proof will follow a localized
version of Viterbo's argument \cite{viterbo2} after proving a simple
lemma in Morse theory. Here we use the constructions of the looping
process of Hamiltonian paths and of Lagrangian suspension together
with the arguments in geometric measure theory in a more systematic
way than Viterbo did in \cite{viterbo2}.

Specialized to the case $U = M$ and stripping the usage of conormals
of general open subsets but restricting to the zero section of the
whole manifold in the Lagrangianization, our proof of the uniqueness
theorem, when restricted to the closed case, is essentially the same
as the one given by Viterbo in \cite{viterbo2} except that we put
his scheme in general perspective of $C^0$-symplectic topology and
$C^0$-Hamiltonian dynamics, and added some more details in the proof
so that we can safely apply them to the current local setting. The
more general case of $L^{(1,\infty)}$-Hamiltonian flows is a subject of
future study.

We will assume that $M$ is closed throughout the paper until section
\ref{sec:general} where we explain how to generalize the main
theorem to  $(M,\omega)$ that is bounded at infinity in the sense of
Gromov \cite{gromov:pseudo}.

I like to thank S. M\"uller and C. Viterbo for useful communications
regarding the proof of Lemma \ref{epsilon} in Appendix.

\bigskip

\noindent{\bf Notations}
\begin{enumerate}
\item $S^1 = \R/\Z$, $S^1(2) = \R / 2 \Z$
\item $\pi: T^*N \to N$, the canonical projection of cotangent bundle.
\item $\widetilde \pi_1: T^*(S^1(2) \times Y) \to T^*S^1(2)$,
$\widetilde \pi_2: T^*(S^1(2) \times Y) \to T^*Y$ the obvious projections.
\item $\pi_1: S^1(2) \times Y \to S^1(2)$, $\pi_2: S^1(2) \times Y \to Y$,
the obvious projections.
\item $\iota_Y: [0,1] \times Y \to T^*Y$, the map defined by $\iota_Y(t,y) = o_y$
where $o_y$ the zero at $y \in T_y^*Y$. Or when $Y \subset (X,\omega)$
is a Lagrangian submanifold, $\iota_Y$ also denotes the corresponding map
$[0,1] \times Y \to X$ induced by the inclusion.
\item $i_Y: Y \to T^*Y,\, X$ the obvious inclusion maps of $Y$.
\end{enumerate}

\section{Preliminaries}
\label{sec:prelim}

\subsection{Conormal of open subsets}
\label{subsec:conormal}

In this section, we recall the definition of the \emph{conormal}
of an open subset and state some results from a
Morse theory, which will be used in the proof of
the locality theorem later.

\begin{defn} Let $U \subset M$ be an open subset with smooth boundary
$\overline U$. We define the \emph{conormal} of $\del U$ by
$$
o_U\coprod\nu^*_-(\partial U)
:= \nu^*\overline{U}
$$
\end{defn}
Here we define $\nu^*_-(\partial U)$ to be the ``negative'' part of
the conormal bundle of $\partial U$, i.e.,
$$
\nu^*_-(\partial U)=\{\alpha\in\nu^*(\partial U)\mid \alpha (\vec n)\le 0,
\vec n\;\hbox{ outward normal to $\partial U$}\}.
$$
In the case $U= (-1,1)\subset\R$, the conormal to $\overline{U} =[-1,1]$
is the one pictured below

\includefigure{kaohfig1}{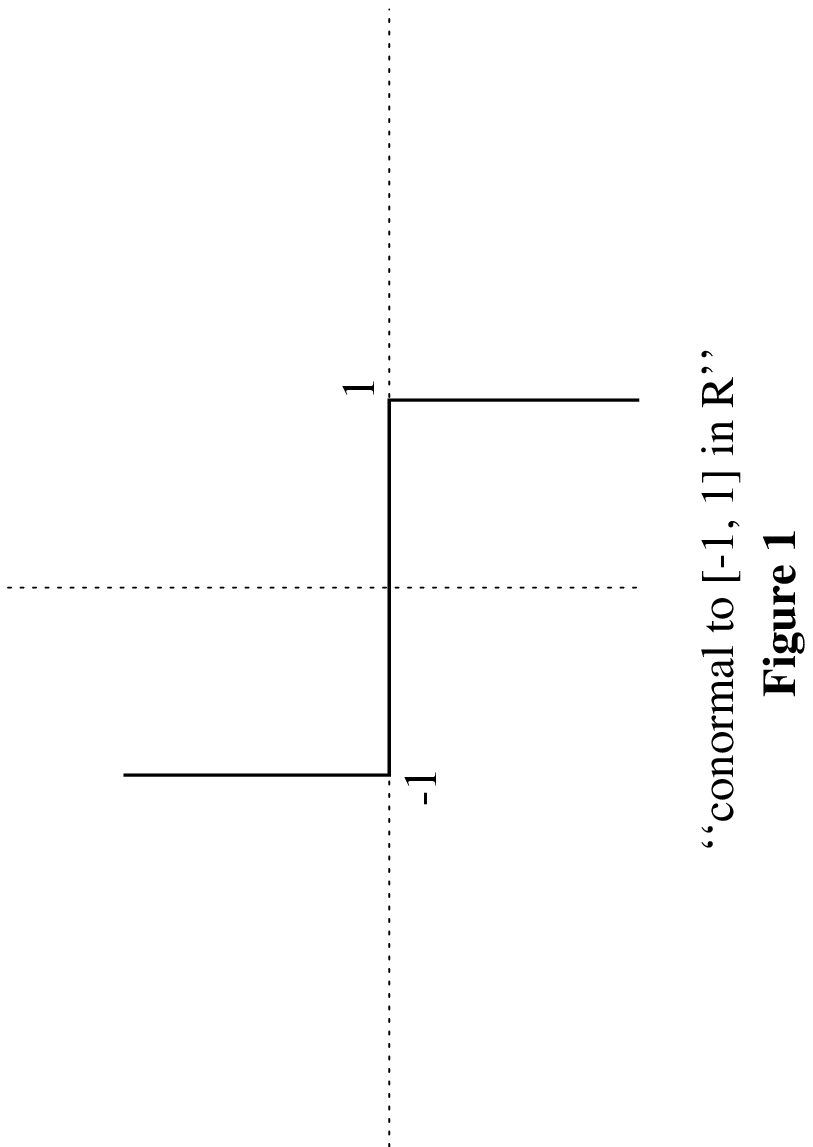}{ht}


The following lemma will be used in the proof of the main theorem.
For the simplicity but without loss of any generality, we assume that
$\overline U$ is compact.

\begin{lem}\label{Morsek} Let $U$ be an open subset of $Y$ with smooth boundary
and let $V$ be an open subset of $Y$ such that
$$
V \subset \overline V \subset U.
$$
Then there exists a smooth Morse function $k : \overline U \to \R$ such that
\begin{enumerate}
\item we have
$$
-dk(\vec n)(y) > 0 \quad \mbox{\rm on }\, \partial U \quad
\mbox{or equivalently }\, \operatorname{Graph}dk \cap \nu^*_-
(\partial U)
$$
\item the function $k : U \to \R$ has no critical
points on $\overline U \setminus \overline V$.
\end{enumerate}
\end{lem}
\begin{proof} We choose any function $k$ constant on $\partial U$
such that (1) holds on $\partial U$. Then perturb the function to
a Morse function away from a neighborhood of $\partial U$.
Finally we compose the function $k$ with
a diffeomorphism of $U$ that is supported away from $\partial U$ and
moves all the critical points into $V$. This finishes the proof.
\end{proof}

\subsection{Continuous Hamiltonian flows and their Hamiltonians}
\label{subsec:topham}

Following \cite{oh:hameo1}, we denote by
$$
\CP^{ham}(Symp(M,\omega),id)
$$
the set of smooth Hamiltonian paths $\lambda: [0,1] \to Symp(M,\omega)$
with $\lambda(0) = id$, and equip it with the
\emph{Hamiltonian topology}. We refer to \cite{oh:hameo1}, \cite{oh:hameo2}
for an extensive study of this topology in the $L^{(1,\infty)}$-context
and \cite{mueller} in the $L^\infty$ or $C^0$-context. In this paper,
we will exclusively use the $C^0$-Hamiltonian topology.

It is equivalent to the metric
topology induced by the metric
$$
d_{ham_\infty}(\lambda,\mu): = \overline d(\lambda,\mu) +
\operatorname{leng}_\infty(\lambda^{-1}\mu)
$$
where $\operatorname{leng}_\infty$ is the $L^\infty$-Hofer length
given by
$$
\operatorname{leng}_\infty(\lambda)
= \|H\|_\infty, \quad \lambda = \phi_H
$$
(Proposition 3.10 \cite{oh:hameo1}), and $\overline d$ is
the $C^0$ metric on $\CP(Homeo(M),id)$.
We consider the \emph{developing map}
$$
\operatorname{Dev}: \CP^{ham}(Symp(M,\omega),id) \to C_m^\infty([0,1] \times M,\R) :
$$
This is defined by the assignment of the normalized generating
Hamiltonian $H$ of $\lambda$, when
$\lambda = \phi_H: t \mapsto \phi_H^t$.
We also consider the inclusion map
$$
\iota_{ham}: \CP^{ham}(Symp(M,\omega),id) \to \CP(Symp(M,\omega),id)
\to \CP(Homeo(M,\omega),id).
$$
Following \cite{oh:hameo1}, we call the product map
$(\iota_{ham},\operatorname{Dev})$ the \emph{unfolding map} and
denote the image thereof by
\be\label{eq:QQ} \CQ :=
\operatorname{Image} (\iota_{ham},\operatorname{Dev}) \subset
\CP(Homeo(M),id) \times C^0_m([0,1] \times M,\R).
\ee
Then both
maps $\operatorname{Dev}$ and $\iota_{ham}$ are Lipschitz with
respect to the metric $d_{ham}$ on $\CP^{ham}(Symp(M,\omega),id)$
by definition and so the unfolding map canonically extends to the
closure $\overline{\CQ}$ in that we have the following continuous
projections
\bea \overline{\iota}_{ham} & : \overline{\CQ} \to
\CP(Homeo(M),id) \label{eq:bariota}\\
\overline{\operatorname{Dev}} & : \overline{\CQ} \to
C^0_m([0,1]\times M,\R). \label{eq:barDev}
\eea

We would like to note that by definition we also have the extension of
the evaluation map $ev_1: \CP^{ham}(Symp(M,\omega),id) \to Symp(M,\omega)
\to Homeo(M)$ to
\be\label{eq:extendev1}
\overline{ev}_1: \operatorname{Image }(\overline\iota_{ham}) \to Homeo(M).
\ee
In this context, Theorem \ref{uniqueness} is
equivalent to saying that the map $\overline{\iota}_{ham}$
is one-one.

Now we recall the notion of \emph{continuous Hamiltonian paths}
studied in \cite{oh:hameo1}, \cite{viterbo2}, \cite{mueller}.

\begin{defn}[Definition 3.18 \cite{oh:hameo1}, \cite{mueller}]
We define the set $\CP^{ham}_\infty(Sympeo(M),id)$ by
$$
\CP^{ham}_\infty(Sympeo(M,\omega),id) : =
\operatorname{Image}(\overline{\iota}_{ham}) \subset
\CP(Homeo(M),id)
$$
and call any element thereof a {\it continuous Hamiltonian path}.
We also define the group of $C^0$-Hamiltonian homeomorphisms,
denoted by $Hameo_\infty(M,\omega)$, to be
$$
Hameo_\infty(M,\omega) := \operatorname{Image}(\overline{ev}_1).
$$
\end{defn}

Similarly to the case of the interval $[0,1]$, we can define a
continuous Hamiltonian path on $[a,b]$ with $b > a$
$$
\lambda : [a,b] \to Hameo(M,\omega)
$$
to be a path such that \be\label{eq:lambda(a)} \lambda \circ
(\lambda(a))^{-1} \in
\CP^{ham}_{\infty;[a,b]}(Sympeo(M,\omega),id) \ee where
$\CP^{ham}_{\infty;[a,b]}(Sympeo(M,\omega),id)$ is defined in the
same way as $\CP^{ham}_\infty(Sympeo(M,\omega),id)$ with $[0,1]$
replaced by $[a,b]$.

We consider smooth functions $\zeta:[0,1] \to [0,1]$ satisfying
\be\label{eq:zeta}
\zeta (t) = \begin{cases}
             0  & \mbox{for }\, 0 \leq t \leq {\e_0}\\
             1  &  \mbox{for }\, 1 - {\e_0} \leq t \leq 1
             \end{cases}
\ee
and
$$
\zeta'(t) \geq 0 \quad \mbox{for all} \quad t \in [0,1],
$$
for some small $\e_0 > 0$, and denote the reparameterized Hamiltonian by $H^\zeta$
defined by
$$
H^{\zeta} (t,x) = \zeta'(t)H(\zeta(t),x)
$$
which generates the Hamiltonian isotopy $ t \mapsto \phi_H^{\zeta(t)}$
in general. Using this reparmetrization, we can always make the
Hamiltonian path \emph{boundary flat} and its Hamiltonian
vanishes near $t = 0, \, 1$.

We will always
assume this boundary-flatness until section \ref{sec:general}, where we will
explain how to treat the general case.

Now we introduce the notion of {\it continuous Hamiltonians}.

\begin{defn}\label{topoham}
We denote
$$
\CH^0_m([0,1] \times M):= \text{Im }(\overline{\operatorname{Dev}})
$$
and call an element in $\CH^0_m([0,1] \times M)$ a
normalized {\it continuous Hamiltonian}.
\end{defn}

By definition, if $H \in
\CH^0_m([0,1]\times M)$, we have a sequence of
smooth Hamiltonian functions $H_i$ such that
\begin{enumerate}
\item $\phi_{H_i}^t$ converges uniformly over $[0,1] \times M$.
\item $H_i$ converges in $C^0_m([0,1] \times M,\R)$ to $H$.
\end{enumerate}

\section{The first half of Theorem \ref{local-uniqueness}}
\label{sec:half}

In this section, we prove a half of Theorem \ref{local-uniqueness}.
This is the localized version of the corresponding theorem from
\cite{viterbo1}, \cite{hofer-zehn}, \cite{oh:hameo1} where $U =M$ is
considered.

\begin{thm}\label{hald}
Let $U \subset M$ be an open subset. Suppose $H_i:[0,1] \times M
\to \R$ is a sequence of smooth normalized Hamiltonian functions
such that both $\phi_{H_i}$ and $H_i$ converge in $C^0$-topology.
Denote by $\lambda$ and $H$ the corresponding $C^0$-limits
respectively. Then if $H(t,x) \equiv c(t)$ on $U$ for a continuous
function $c: [0,1] \to \R$ depending only on $t$, then $\lambda
\equiv id$ on $[0,1] \times U$.
\end{thm}
\begin{proof} We first note that $H(t,x) \equiv c(t)$ on $U$
is equivalent to $\|H\|_{\infty,U} = 0$.

It is enough to prove that for any given
$\delta>0$, we have
\be\label{eq:0onUdelta}
\lambda(t,x) = x \quad \mbox{for all $x \in U_\delta$}
\ee
where $U_\delta$ is the subset
$$
U_\delta = \{x \in U \mid d(x,\partial U) > \delta\}.
$$
Suppose (\ref{eq:0onUdelta}) does not hold and so there exists $t_0
\in [0,1]$ and $x_0 \in U_\delta$ such that $\lambda(t_0,x_0) \neq
x_0$. Then there exists $\e_0 > 0$ such that the closed ball
$B_{\e_0}(x_0)$ satisfies
$$
\lambda_{t_0}(B_{\e_0}(x_0)) \cap B_{\e_0}(x_0) = \emptyset.
$$
Furthermore, since $\lambda$ is continuous and $\lambda_0 = id$
there exists $t_1 > 0$ such that
$$
t_1: = \max\{t \in [0,1] \mid \lambda_s(B_{\e_0}(x_0)) \subset \overline U
\, \, \mbox{for all $0 \leq s \leq t_1$}\}.
$$
Then we can choose another smaller ball $B_{\e_1}(x_1) \subset
B_{\e_0}(x_0)$ such that
\be\label{eq:Be1x1}
\lambda_{t_1}(B_{\e_1}(x_1)) \cap B_{\e_1}(x_1) = \emptyset
\ee
and
\be\label{eq:Be1x1inU}
\lambda_t(B_{\e_1}(x_1)) \subset U
\ee
for all $0 \leq t \leq t_1$ in addition.

By the uniform convergence of $\phi_{H_i}$, we still have
\be\label{eq:Hidisjoint}
\phi_{H_i}^t(B_{\e_1}(x_1)) \cap B_{\e_1}(x_1) = \emptyset
\ee
for all $0 \leq t \leq t_1$ if we choose a sufficiently large $i$.
Choose an open subset $V \subset  U$ such that $\overline V \subset U$ and
$$
V \supset \bigcup_{0 \leq t \leq t_1} \phi_{H_i}^t(B_{\e_1}(x_1)).
$$
We multiply to $H_i$ a cut-off function $\chi: Y \to \R$ satisfying
$$
\chi(x) = \begin{cases} 1 & \quad x \in \bigcup_{0 \leq t \leq t_1}
\phi_{H_i}^t(B_{\e_1}(x_1)) \\
0 & \quad x \in Y\setminus V
\end{cases}
$$
and define a new Hamiltonian $F_i = \chi H_i$. Then we still have
\be\label{eq:Hidisjoint} \phi_{F_i}^t(B_{\e_1}(x_1)) \cap
B_{\e_1}(x_1) = \emptyset \ee for all $0 \leq t \leq t_1$ because
the flow of $F_i$ and that of $H_i$ coincide on $\bigcup_{0 \leq t
\leq t_1} \phi_{H_i}^t(B_{\e_1}(x_1))$ since  $\chi \equiv 1$
thereon. Therefore we obtain
\be\label{eq:F1t1}
\|F_i^{t_1}\|_\infty \geq e(B_{\e_1}(x_1)) > 0
\ee
where $F_i^{t_1}$ is the reparameterized $F_i$ given by
$$
F_i^{t_1}(t,x) = F_i(t_1t,x)
$$
and $e(B_{\e_1}(x_1))$ is the displacement energy of
$B_{\e_1}(x_1)$. (See \cite{hofer2}, \cite{lal-mc1} for its
definition.) The positivity of $e(B_{\e_1}(x_1))$ follows from the
energy-capacity inequality from \cite{lal-mc1}. On the other hand,
obviously we have
$$
F_i^{t_1} = \chi H_i^{t_1} \to \chi H^{t_1}
$$
in $C^0$-topology and hence (\ref{eq:F1t1}) implies
\be\label{eq:chiH>0}
\|\chi H^{t_1}\|_\infty > \frac{1}{2}e(B_{\e_1}(x_1)) > 0
\ee
for all sufficiently large $i$'s.
However we have
\bea
\label{eq:chiH=0}
\|\chi H^{t_1}\|_\infty & = & \max_{t \in [0,1]}
\left(\max_{x \in M}\osc (\chi H(t_1t,x))\right) \nonumber\\
& = & \max_{t \in [0,t_1]}
\left(\max_{x \in M}\osc (\chi H(t,x))\right) \nonumber \\
& \leq & \max_{t\in [0,1]} \left(\max_{x \in U}\osc H^{t_1}(t,x)\right) =
\|H\|_{\infty,U} = 0 \eea where the inequality holds since
$\operatorname{supp} \chi  \subset U$ for all $t \in [0,1]$ for a
sufficiently $i$ and $0 \leq \chi(t) \leq 1$.

Obviously (\ref{eq:chiH>0})
and (\ref{eq:chiH=0}) contradict to each other.
This finishes the proof.
\end{proof}

For the rest of the  paper, we will provide the other half of
the proof of Theorem \ref{local-uniqueness}.

\section{Continuous exact Lagrangian isotopy}
\label{sec:uniqueness}

In this section, we denote a general symplectic manifold
by $(X,\omega)$, instead of using the letter $M$.
We also denote a closed manifold  by $Y$. When we are given a
Lagrangian embedding of $Y$ into $(X,\omega)$, we often do not
distinguish $Y$ from the image in $X$ of the given embedding.

Let $\psi: [0,1] \times Y \to X$ be an isotopy of embeddings
$\psi_t:Y \to X$. The pull-back $\psi^*\omega$ is decomposed to
\be\label{eq:psiomega}
\psi^*\omega = dt \wedge \alpha + \beta
\ee
where $\alpha$ and $\beta$ are forms on $Y$ satisfying
$$
\frac{\del}{\del t}\rfloor \alpha = 0= \frac{\del}{\del t}\rfloor \beta.
$$
The following definition is given in \cite{gromov:pseudo}.

\begin{defn}\label{exact-lag} An isotopy $\psi: [0,1] \times Y \to X$
is called a Lagrangian isotopy if $\beta = 0$
and $i_t^*\alpha$ is closed (respectively, an \emph{exact Lagrangian isotopy})
if in the above decomposition $i_t^*\alpha$ is closed
(respectively, exact) for all $t \in [0,1]$.
\end{defn}
We remark that in this definition, the condition $\beta = 0$
automatically implies closedness of $i_t^*\alpha$. Therefore the latter
condition is redundant for the definition of Lagrangian isotopy.
We however put it there separately to make comparison with
the case of exact Lagrangian isotopy.

The following lemma is well-known and is an easy consequence of the
definition whose proof we omit.

\begin{lem}\label{psih}
If $\psi$ is an exact Lagrangian isotopy, then there exists a
function $h:[0,1] \times Y \to \R$ such that $\psi^*\omega =
dt \wedge dh$ i.e.,
$$
\psi^*\omega + dt \wedge (-dh) = 0.
$$
Furthermore if $h'$ is another such function, we
have $h' - h$ is a function depending only on $t$, but not on $Y$.
\end{lem}
We call any such $h$ a {\it hamiltonian} of $\psi$. Here we like to
alert readers that we use the lower case `h' to call the
function $h$ the `hamiltonian' associated to $\psi$ to distinguish it from
the Hamiltonian functions defined on the ambient symplectic
manifold $X$.

According to this lemma, for a given exact Lagrangian isotopy $\psi$
its associated hamiltonian is not unique. As we will see later, the
suspension of an exact Lagrangian isotopy depends not just on
the isotopy $\psi$ but also on its hamiltonian. Because of this,
we will regard the pair
$$
(h,\psi): [0,1] \times Y \to \R \times X
$$
as an exact Lagrangian isotopy. With a suitable normalization on
the hamiltonians associated to $\psi$,
the choice of $h$ can be made unique in the smooth category.

The following extension lemma is also well-known whose proof
is omitted.

\begin{lem}\label{exact=ham}
For any exact Lagrangian isotopy $\psi: [0,1] \times Y \to X$ of
Lagrangian embeddings, there exists a Hamiltonian $H: [0,1] \times X \to \R$
such that
\be\label{eq:exact=ham} \psi(t,x) = \phi_H^t(\psi(0,x))
= \phi_H^t(\psi_0(x)). \ee
\end{lem}

When the exact Lagrangian isotopy $(h,\psi)$ is given,
we can always adjust $H$ away from the support of the isotopy
so that $H$ satisfies the normalization condition
$$
\int_X h_t \, \Omega_\omega = 0
$$
where $\Omega_\omega$ is the Liouville volume form.

And if the isotopy is boundary flat, then $h$ so that
$H$ can be made boundary flat. Namely, we may assume that
there exists $\e_0 >0$ such that
\be\label{eq:H=0}
H\equiv 0 \quad \mbox{for } \, t \in [0,\e_0] \cup [1-\e_0,1].
\ee

\begin{defn}\label{reversal} Let $Y \subset X$ be a Lagrangian
submanifold and $(h,\psi)$ an exact Lagrangian isotopy.
The \emph{time-reversal} of $(h,\psi)$ is
the pair $(\widetilde h, \widetilde \psi)$ defined by
$$
\widetilde \psi(t,x) = \psi(1-t,x), \quad \widetilde h(t,x) =
-h(1-t,x).
$$
\end{defn}

Obviously for a boundary flat $(h,\psi)$, the concatenated isotopy
$$
(h,\psi) * (\widetilde h, \widetilde \psi)
$$
defines a \emph{smooth} closed embedding $S^1(2) \times Y
\to \R \times S^1(2) \times X$, where we represent $S^1$ by $\R/2\Z$.
We call this concatenated isotopy an \emph{odd double} of the isotopy
$(h,\psi)$ and denote it by $(h^{od},\psi^{od})$.

Now we are ready to give the definition of the \emph{odd-double suspension} of
an exact Lagrangian isotopy $(h,\psi)$.

\begin{defn}\label{odd-double} The odd-double suspension of an exact Lagrangian isotopy
$(h,\psi): [0,1] \times Y \to \R \times X$ is defined by
the embedding
$$
\iota_{(h,\psi)} : S^1(2) \times Y \to T^*S^1(2) \times X
$$
\be\label{eq:iotahpsi}
\iota_{(h,\psi)}(\theta,x) =
\begin{cases} (\theta, - h(\theta,x),\psi(\theta,x)) \quad
& \mbox{for $\theta \in [0,1]$} \\
(\theta,  h(2-\theta,x), \psi(2-\theta,x)) \quad
& \mbox{for $\theta \in [1,2]$}.
\end{cases}
\ee We also denote the
obvious double $h * \widetilde h$ of the hamiltonian $h$ by
$h^{od}$.
\end{defn}

Now we state a proposition which is an essential ingredient in the
proof of Theorem \ref{local-uniqueness}.
We provide its proof in the appendix.
(See also Theorem 6.1.B \cite{polterov}.)

\begin{prop}\label{isotopeinU} Let $(h,\psi)$ be an exact Lagrangian isotopy.
Suppose that the image of $\psi$ is contained in
$U \subset X$. Equip $T^*S^1(2) \times X$ with the symplectic form
$$
d\theta \wedge db + \omega_X
$$
where $(\theta,b)$ is the canonical coordinates of $T^*S^1(2)$.
Then the Lagrangian embedding $\iota_{(h,\psi)}: S^1(2)\times
Y \to T^*S^1(2) \times X$ is isotopic to the embedding
$o_{S^1(2)} \times i_Y : S^1(2) \times Y \subset T^*S^1(2)
\times X$ by a Hamiltonian flow supported in $T^*S^1(2) \times U$.
\end{prop}

We assume that for
any $s \in [0,1]$, $h(s, \cdot)$ assumes the value zero so that
\be
\min_{x} h_s \leq 0 \leq \max_x h_s.
\ee
This can be done, for example, by normalizing $h_s$ so that
\be\label{eq:average0}
\frac{\max_x h_s + \min_x h_s}{2} = 0.
\ee
In particular, this normalization makes
\be
\min_{x}h_s < 0 < \max_x h_s
\ee
\emph{unless $h_s \equiv 0$.} We call any such hamiltonian $h$
satisfying (\ref{eq:average0}) \emph{average normalized}.

We denote by $L_s = \phi_s(Y)$ the family of Lagrangian submanifolds
associated to the isotopy. Without loss of generality, we will assume that
$\psi$ is boundary flat. We introduce the following definition

\begin{defn}\label{topoexact}
A pair $(h,\psi)$ such that $\psi:[0,1] \times Y \to X$ is a
continuous map and $h:[0,1] \times Y \to \R$ an $C^0$-function is
called a {\it continuous exact Lagrangian isotopy} if there exists
a sequence of smooth exact Lagrangian isotopies $(h_i,\psi_i)$
such that $\psi_i \to \psi$ in the $C^0$-topology and $h_i \to h$
in the $C^0$-topology. We call any such sequence $(h_i,\psi_i)$
a \emph{smoothing sequence} of $(h,\psi)$.
\end{defn}

Similarly as the definition of $\CH^0_m([0,1]\times M)$, we
give the following Lagrangian counterpart thereof.

\begin{defn}\label{topoham} We define the subset
$$
h^0_{av}([0,1] \times Y) \subset C^0([0,1] \times L)
$$
to be the set of the above $C^0$-limits of normalized sequences of
smooth exact Lagrangian isotopies $(h_i,\psi_i)$. We call any
element $h \in h^0_{av}([0,1] \times Y)$ a
(time-dependent) $C^0$-hamiltonian function on $Y$.
\end{defn}

We define the norms
$$
\|h\|_\infty = \max_{t \in [0,1]}\left(\max_{x\in Y} h_t - \min_{x\in Y}
h_t\right)
$$
and for a given open subset $U$, we define
$$
\|h\|_{\infty,U} = \max_{t \in [0,1]} \left(\max_{x\in U} h_t -
\min_{x\in U} h_t\right)
$$
for an average-normalized hamiltonian $h$.

The following theorem is the Lagrangian version of Theorem
\ref{local-uniqueness}, which will be proved in the rest of the paper.

\begin{thm}\label{lag-localuniqueness}
Let $(X,\omega)$ be a symplectic manifold and $Y \subset X$ be a
closed Lagrangian submanifold. We denote by
$$
\iota_Y: [0,1] \times Y \to [0,1] \times X
$$
the obvious inclusion map. Let $(h_i,\psi_i): [0,1] \times Y \to \R \times X$ be a
sequence of exact Lagrangian isotopies that are boundary flat.
Let $U \subset Y$ be an open subset with smooth boundary $\del U$.

Suppose that $(h_i,\psi_i)$ converges to $(h,\psi)$ in $C^0$-topology
on $Y$. Then if $\psi \equiv \iota_Y$
on $[0,1] \times U$, we have $h(t,y) \equiv c(t)$ for $y \in U$
where $c:[0,1] \to \R$ is a continuous function depending only on
$t$.
\end{thm}

\section{Lagrangian intersections and analysis of Calabi currents}
\label{sec:analysis}

We assume that $Y$ is oriented throughout this section. The proof for
the unoriented case will be similar if one uses \emph{odd currents}
\cite{derham} whose details we omit because it will not be used
in the proofs of the main theorems in this paper.

We first state the following lemma whose proof will be given in section
\ref{sec:notclosed} in the more general $L^{(1,\infty)}$-context.

\begin{lem}\label{notclosed} Let $\alpha = h^{od} d\theta$
be the Calabi current associated to the continuous exact Lagrangian
isotopy $(h,\psi)$ as defined above. Suppose $\|h\|_{\infty,U} \neq
0$. Then the current $\alpha$ on $S^1(2) \times Y$ is not closed on
$U$.
\end{lem}

The following definition will be useful for our purpose.

\begin{defn}\label{phi} A diffeomorphism $\phi: S^1(2) \times Y \to S^1(2) \times Y$
is called \emph{fiberwise-like} with respect to the projection
$\pi_2: S^1(2) \times Y \to Y$, if it has the form
$$
\phi(\theta,y) = (\theta, Y_\theta(y))
$$
for a 2-periodic family of diffeomorphisms $Y_\theta: Y \to Y$
with $Y_0 = id_Y$. We define its support by
$$
\operatorname{supp}\phi = \bigcup_{\theta\in [0,2]}\operatorname{supp}
Y_\theta.
$$
\end{defn}

Obviously the space $C^0(S^1(2) \times Y)$ is invariant under a
full diffeomorphism group of $S^1(2) \times Y$.

For each fiberwise-like diffeomorphism $\phi: S^1(2) \times Y \to
S^1(2) \times Y$, we consider the push-forward $\phi_*\alpha$ as a
current and decompose it into
\be\label{eq:alphaphi} \phi_*\alpha
= \beta_\phi + (\phi_*\alpha - \beta_\phi)
\ee
where $\beta_\phi$ satisfies $\frac{\del}{\del \theta} \rfloor
\beta_\phi = 0$. Note that $\phi_*\alpha - \beta_\phi$ can be
written as
$$
\phi_*\alpha - \beta_\phi = \left(\frac{\del}{\del \theta} \rfloor
\phi_*\alpha\right) \wedge d \theta.
$$
The following can be proved in a straightforward calculation and so
its proof is omitted. We however point out that $\phi_*\alpha$
defines a well-defined $n$-current on $S^1(2) \times Y$.

\begin{lem}\label{phi-invariance} Let $\phi$ be a fiberwise-like
diffeomorphism. Under the decomposition (\ref{eq:alphaphi}), we
have
$$
\phi_*\alpha - \beta_\phi = h_\phi^{od}\, d \theta
$$
where $h_\phi^{od}$ is an $C^0$-function which is given by the
formula \be\label{eq:hphiod} h_\phi^{od}(\theta,y) : =
h^{od}(\phi^{-1}(\theta,y)) = h^{od}(\theta,Y_\theta^{-1}(y)). \ee
\end{lem}

Using this preparation, we state the following proposition whose
proof will be given in section \ref{sec:transport} in the more general $L^{(1,\infty)}$-context.
We also refer the readers to Lemma 2.5 \cite{viterbo2} for the
corresponding global statement.

\begin{prop}\label{alphaphi}
Let $U \subset Y$ be an open subset.
Let $\alpha = h^{od}\, d\theta$ be the above mentioned current on
$S^1(2) \times Y$. Suppose that $\alpha$ is not closed on $U$. Then there
exists a fiberwise-like diffeomorphism $\phi: S^1(2) \times Y \to
S^1(2) \times Y$ supported on $S^1(2) \times U$
such that we have some $y \in U$ such that
$$
\int_{\pi_2^{-1}(y)} h_\phi^{od}\, d\theta \neq 0.
$$
Furthermore $\phi$ can be chosen arbitrarily $C^\infty$-close
to the identity map and $\operatorname{supp} \phi \subset U$.
\end{prop}

With Lemma \ref{notclosed} and Proposition \ref{alphaphi} assumed,
we now start the proof of Theorem \ref{lag-localuniqueness}.

For this purpose, we closely analyze the function $P_{h_\phi^{od}}:
Y \to \R$ on $U$ which is defined by \be\label{eq:Pphi}
P_{h_\phi^{od}}(y) = \int_{\pi_2^{-1}(y)} h_\phi^{od}\, d\theta \ee
for $y \in Y$. This function is obviously continuous on $Y$
for a $C^0$-function $h$ and for a diffeomorphism $\phi$. Here we
like to note that due to the doubling process, the function $P$
given by
$$
P_{h^{od}}(y) = \int_{\pi_2^{-1}(y)} h^{od}\, d\theta
$$
itself identically vanish, which prompts one to deform it by a
diffeomorphism or to deform the fibrations as put in
\cite{viterbo2}.

The following will be used in our proof in an essential way.

\begin{lem}\label{Pphi}
Let $\alpha = h_\phi^{od}d\theta$ with $h \in C^0$ and
$\|h_\phi^{od}\|_{\infty,U} \neq 0$. Then $P_{h_\phi^{od}}: Y \to
\R$ is a non-zero continuous function such that
\begin{enumerate}
\item
the set
$$
U_\phi: =\{y \in \overline U \mid P_{h_\phi^{od}}(y) \neq 0 \}
$$
is a non-empty open subset of $\overline U$.
\item
$P_{h_\phi^{od}} \equiv 0$ on $\overline U \setminus
\operatorname{supp}\phi \supset \partial U$ where $\phi(\theta,y) =
(\theta,Y_\theta(y))$.
\end{enumerate}
\end{lem}
\begin{proof}
We start with the formula
$$
P_{h_\phi^{od}}(y) = \int_{S^1(2)} h^{od}(\phi^{-1}(\theta,y)) \,
d\theta = \int_{S^1(2)} h^{od}(\theta,Y_\theta^{-1}(y)) \, d\theta
$$
which follows from (\ref{eq:hphiod}). Obviously $P_{h_\phi^{od}}$ is
a continuous function on $Y$ which satisfies $P_{h_\phi^{od}} \not
\equiv 0$ by Proposition \ref{alphaphi}.

Next, we recall that
$$
\int_0^2 h^{od}(\theta,y) \, d\theta = \int_0^1 h(\theta,y)\, d\theta
+ \int_1^2(-h(2-\theta,y))\,d\theta \equiv 0
$$
for all $y \in Y$. Therefore we derive
$$
P_{h_\phi^{od}}(y) = \int_0^2
h^{od}(\theta,Y_\theta^{-1}(y)) \, d\theta = 0
$$
whenever $y \in U \setminus \operatorname{supp}\phi$, i.e.,
$Y_\theta^{-1}(y) = y$ for all $\theta \in [0,2]$.


This finishes the proof.
\end{proof}

The following variant of H\"ormander's construction
(see (2.5.4) \cite{hormander})
will be useful for our discussion coming henceforth.

\begin{defn}\label{N1-N2}
Consider the cotangent bundle $T^*N$ of a smooth manifold.
Let $\psi_i: N_i \to T^*N$ be two Lagrangian embeddings of manifolds
$N_i$, $i = 1,\,2$.
We introduce the \emph{difference set} denoted by $\psi_1 \ominus \psi_2$.
The set $\psi_1 \ominus \psi_2$ is defined by
\be\label{eq:N1-N2}
\psi_1 \ominus \psi_2 =\{ (q,p) \in T^*N \mid
(q,p) = \psi_1(n_1) - \phi_2(n_2), \, q \in N, \, \psi_i(n_i)
\in N_i \cap T_q^*N\}
\ee
\end{defn}

Now consider the symplectic diffeomorphism
$$
T^*\phi^{-1}: T^*(S^1(2) \times Y) \to T^*(S^1(2) \times Y)
$$
given by
$$
T^*\phi^{-1}(\alpha_{(\theta,y)}) =
(T_{\phi(\theta,y)}\phi^{-1})^*(\alpha_{(\theta,y)})
\in T^*_{\phi(\theta,y)}(S^1(2) \times Y).
$$
The composition
$$
\iota_{(h_i,\psi_i)}^\phi : = T^*\phi^{-1} \circ \iota_{(h_i,\psi_i)} \circ \phi^{-1}:
S^1(2) \times Y \to T^*(S^1(2) \times Y)
$$
defines a Lagrangian embedding. Using the splitting
$T^*(S^1(2) \times Y) = T^*S^1 \times T^*Y$ we write
\be\label{eq:iotaetai} \iota_{(h_i,\psi_i)}^\phi(\theta,y) =
(\alpha_i^\phi(\theta,y),\beta_i^\phi(\theta,y)) \ee as an element
in $T^*_{(\theta',y')}(S^1(2) \times Y)$ where
$$
(\theta',y') = \pi(\iota_{(h_i,\psi_i)}^\phi(\theta,y)).
$$
Here we note that $\theta' = \theta$ as $\phi$ has the form
$\phi(\theta,y) = (\theta,Y_\theta(y))$. We denote by $\widetilde
\pi_i, \, i = 1,\,2$ the projections from $T^*(S^1(2) \times Y)$ to
$T^*S^1$ and $T^*Y$ respectively.

We note the identity
$$
\pi\circ \iota_{(h_i,\psi_i)}^\phi = \pi \circ \psi_i.
$$
since $\psi_i$ converges to the zero section map $\iota_Y: S^1(2)
\times Y \to T^*(S^1(2) \times Y)$ in $C^0$ topology, the map
$$
\pi\circ \iota^\phi_{(h_i,\psi_i)}: S^1(2) \times Y \to S^1(2)
\times Y
$$
can be made arbitrarily close to the identity uniformly over $i$.

Similarly we define a continuous embedding $\iota_{(h,\iota_Y)}^\phi:
S^1(2) \times Y \to T^*Y$ by
\be\label{eq:iotahid}
\iota_{(h,\iota_Y)}^\phi(\theta,y) = (h\circ \phi^{-1},id).
\ee
The following proposition shows that $\iota_{(h,\iota_Y)}^\phi$ is again
a continuous \emph{Lagrangian} embedding.

\begin{prop}\label{oS1Y} The smooth Lagrangian embedding
$\iota_{(h_i,\psi_i)}^\phi$ converges to the continuous embedding
$\iota_{(h,\iota_Y)}^\phi$ in $C^0$-topology.
\end{prop}
\begin{proof} For this, we compute the values of $\alpha_i^\phi(\theta,y)$ as an element
in $T^*_{(\theta',y')}(S^1(2) \times Y)$ explicitly.

\begin{lem}\label{hiphi}  We have
\be\label{eq:alphaiphi} \alpha_i^\phi(\theta,y) =
h_i(\phi^{-1}(\theta,y))\, d\theta + \widetilde \pi_1\circ
T^*\phi^{-1}\circ \psi_i(\phi^{-1}(\theta,y)).
\ee
\end{lem}
\begin{proof} We compute
\beastar \iota_{(h_i,\psi_i)}^{\phi}(\theta,y)
\left(\frac{\del}{\del \theta}\right) & = & T^*\phi^{-1}\circ
\iota_{(h_i,\psi_i)}(\phi^{-1}(\theta,y))
\left(\frac{\del}{\del \theta}\right) \\
& = & \iota_{(h_i,\psi_i)}(\phi^{-1}(\theta,y))
\left(T\phi^{-1}\left(\frac{\del}{\del \theta}\right)\right)\\
& = & h_i^{od}(\phi^{-1}(\theta,y)) d\theta \left(T\phi^{-1}
\left(\frac{\del}{\del \theta}\right)\right) + \psi_i^{od}((\phi^{-1}(\theta,y))
\left(T\phi^{-1}\left(\frac{\del}{\del \theta}\right)\right)
\eeastar
where we used $\alpha_i = h_i^{od}\, d\theta$.

Using the identity $\phi^{-1}(\theta,y) = (\theta, Y_\theta^{-1}(\theta,y))$,
the first term becomes
$$
h_i^{od}(\phi^{-1}(\theta,y)) d\theta \left(\frac{\del}{\del
\theta}\right) = h_i^{od}(\phi^{-1}(\theta,y))
$$
which gives rise to the first term in (\ref{eq:alphaiphi}).

The second term here gives rise to the second one in
(\ref{eq:alphaiphi}). Hence the proof.
\end{proof}

Obviously the first term of (\ref{eq:alphaiphi}) converges to
$h_i\circ \phi^{-1} \, d\theta$. We now estimate the second term.

\begin{lem}\label{4piC}
Denote by $\iota_Y: S^1(2) \times Y \to X$ the obvious
inclusion-induced map given by $(\theta,y) \mapsto y$. Let
$(h_i,\psi_i) \to (h,\iota_Y)$ in the hamiltonian topology and
denote
$$
m_i(\phi): = \max_{(\theta,y)}\, \operatorname{dist}
\left(\widetilde \pi_1\circ T^*\phi^{-1}(\psi_i(\phi^{-1}(\theta,y)),o_{S^1}
\right).
$$
Then there exits a constant $C>0$ independent of $\phi$ and $i$ such
that \be\label{eq:4piC} m_i(\phi) \leq C d_{C^0}(\psi_i,\iota_Y)
|\phi^{-1}|_{C^1} \ee for all $i$. In particular, $m_i(\phi) \to 0$
as $i \to \infty$.
\end{lem}
\begin{proof}
We have
\beastar m_i(\phi) & \leq &
\max_{(\theta,y)}\left|\widetilde \pi_1\circ T^*\phi^{-1}\circ
\psi_i(\phi^{-1}
(\theta,y))\left(\frac{\del}{\del \theta}\right) \right|\\
& = & \max_{(\theta,y)}\left|\psi_i(\phi^{-1}(\theta,y))
\left(T\phi^{-1}\left(\frac{\del}{\del \theta}\right)\right)\right| \\
& \leq & C \max_{(\theta,y)}|\psi_i(\theta,y)||\phi^{-1}|_{C^1}
\eeastar
for a universal constant $C>0$ where $|\psi_i(\phi^{-1}(\theta,y))|$
is the norm taken as an element of $T^*_{y'}(Y)$ with
$$
y' = \pi_2\phi(\pi(\psi_i(\phi^{-1}(\theta,y))))
= Y_\theta(\pi(\psi_i(\theta,Y_\theta^{-1}(y)))).
$$
 Here we recall that we identify the exact Lagrangian
isotopies $\psi_i: S^1(2) \times Y \to X$ with the corresponding
isotopies in the Darboux chart $T^*Y$ of $Y \subset X$. Under this
identification, we can identify
$$
\max_{(\theta,y)}|\psi_i(\theta,y)|
$$
with $d_{C^0}(\psi_i,\iota_Y)$. This finishes the proof of the first
statement. The last statement follows from (\ref{eq:4piC}) by the
$C^0$-convergence $\psi_i \to \iota_Y$ that we assumed.

This finishes the proof of the lemma.
\end{proof}

Combining Lemma \ref{hiphi} and Lemma \ref{4piC}, we have finished the
proof of Proposition \ref{oS1Y}.
\end{proof}

%

\section{Proof of Theorem \ref{lag-localuniqueness}}
\label{sec:wrap-up}

In this section, we will wrap-up the proof of Theorem
\ref{lag-localuniqueness} by producing a contradiction to Theorem
\ref{kasturi-oh}, \emph{if we assume $\|h\|_{\infty,U} \neq 0$}.


For the simplicity of notation, we denote
$$
P_\phi = P_{h_\phi^{od}}.
$$
%

The following is the local version of the key ingredient
used in Viterbo's scheme \cite{viterbo2}.

\begin{prop}\label{continuoush}
Suppose $\|h\|_{\infty,U} \neq 0$. Let $\phi$ be a fiberwise-like
diffeomorphism as in Proposition \ref{Pphi}. Then there exists
 a smooth function $f: S^1(2) \times Y \to \R$ such that
$$
\iota_{(h,\iota_Y)}^\phi \ominus df
$$
does not intersect $\nu^*(S^1(2) \times \overline U)$.
\end{prop}

Once we have proved this proposition,
an immediate corollary of Proposition {\ref{oS1Y} will be the following
intersection result.

\begin{cor}\label{non-intersect}
Suppose $\|h\|_{\infty,U} \neq 0$ and let $f$ be the smooth function
as in Proposition \ref{continuoush}. Then there exists an $N \in \N$
for which the smooth Lagrangian embedding $\iota_{(h_N,\psi_N)}^\phi
\ominus df$ cannot intersect $\nu^*(S^1(2) \times \overline U)$.
\end{cor}
\begin{proof} Since taking the difference $\ominus df$ is a
continuous operation in $C^0$-topology, Proposition
\ref{oS1Y} implies that $\iota_{(h_i,\psi_i)}^\phi \ominus df$
converges to $\iota_{(h,\iota_Y)}^\phi \ominus df$ in $C^0$-topology.
Therefore Proposition \ref{continuoush} and compactness of
$S^1(2) \times Y$ imply that $\iota_{(h_N,\psi_N)}^\phi \ominus df$
do not intersect $\nu^*(S^1(2) \times \overline U$ if we choose
a sufficiently large $N$.
\end{proof}

With this corollary, we can wrap-up the proof of Theorem \ref{lag-localuniqueness}.

\begin{proof}[Proof of Theorem \ref{lag-localuniqueness}]
\emph{If we assume $\|h\|_{\infty,U} \neq 0$}, Proposition
\ref{continuoush} holds and so Corollary \ref{non-intersect} would
give rise to a contradiction to Theorem \ref{kasturi-oh} because
$\iota_{(h_N,\psi_N)}^\phi \ominus df$ is Hamiltonian isotopic to
the zero section : the latter
follows since $\phi$ is isotopic to the identity on $S^1(2) \times
Y$. and $\iota_{(h_N,\psi_N)}$ is Hamiltonian isotopic to the zero
section. This will then finish the proof of Theorem
\ref{lag-localuniqueness}.
\end{proof}

The rest of the section will be occupied with the proof of
Proposition \ref{continuoush}.

%

\begin{proof}[Proof of Proposition \ref{continuoush}]

Motivated by Viterbo's scheme used in the proof of Proposition 2.1
\cite{viterbo2}, we look for a smooth function $f = f(\theta,y)$
such that
\be\label{eq:difference} \left(\iota_{(h,\iota_Y)}^\phi \ominus
df\right) \cap \nu^*(S^1(2) \times \phi(\overline U)) =
\left(\iota_{(h,\iota_Y)}^\phi \ominus df\right) \cap
\nu^*(S^1(2) \times \overline U) = \emptyset
\ee
where the set
$\iota_{(h,\iota_Y)}^\phi \ominus df$ is the difference set
introduced in Definition \ref{N1-N2}. And we also use the fact that
$\operatorname{supp}\phi \subset U$ and so
$$
\phi(\overline U) = \overline U
$$
for the first identity.

We first solve
\be\label{eq:dfdt}
\begin{cases}
h^{od}_\phi(\theta,y) - \frac{\del \widetilde f}{\del
\theta}(\theta,y)
= P_\phi(y) \\
\widetilde f(0,y) = 0
\end{cases}
\ee for some $\widetilde f$. This is uniquely solvable for given $y
\in Y$ for a smooth function $\widetilde f(\cdot, y)$ because we
have \be\label{eq:inthe} \int_0^1 (h^{od}_\phi(\theta,y) -
P_\phi(y))) \, d\theta = 0 \ee by the definition of $P_\phi$. We
denote by $\widetilde f = \widetilde f(\theta,y)$ such a solution.
However $\widetilde f$ is a priori only continuous and hence we will
suitably perturb $P_\phi$ so that the resulting solution becomes
smooth. In this regard, the following lemma will be important. This
is what Viterbo briefly mentioned in the proof of Lemma 2.4
\cite{viterbo2}. Because the lemma is not totally obvious, we will
give its proof in the appendix for completeness' sake. To make
comparison of our proof with those in \cite{viterbo2} easier, we will use
the same letter $\e$ for the perturbed function as in \cite{viterbo2}.

\begin{lem}\label{epsilon} There exists a continuous function $\e: S^1(2) \times Y
\to \R$ such that
\begin{enumerate}
\item For all $y \in Y$, $\int_{S^1(2)} (h_\phi^{od}(t,y) - \e(t,y)) \, dt = 0$,
\item $h_\phi^{od}-\e$ is smooth.
\item We have the inequality
\be\label{eq:<1/30} \|P_\phi - \e \|_{C^0} <
\frac{1}{30}\|P_\phi\|_{C^0}. \ee
\end{enumerate}
\end{lem}

Now we consider the equation
\be\label{eq:dfdt}
\begin{cases}
h^{od}_\phi(\theta,y) - \frac{\del g}{\del
\theta}(\theta,y)
= \e(\theta,y) \\
 g(0,y) = 0
\end{cases}
\ee
instead. Then the solution
$$
g(\theta,y) = \int_0^\theta (h_\phi^{od}(t,y) - \e(t,y))\, dt
$$
becomes $2$-periodic and so well-defined on
$S^1(2) \times Y$ and becomes smooth by the choice of $\e$ in
Lemma \ref{epsilon}.

We now choose $\delta $ so small independently of $i$, (e.g., the choice
$$
\delta = \frac{1}{10} \|P_\phi\|_{C^0}
$$
will do the purpose) that for the set
$$
Z_{P_\phi}^\delta = \{(\theta,y) \mid |P_\phi(y)| \leq \delta \}
$$
$Y \setminus \pi(Z_{P_\phi}^\delta)$ is still a non-empty open
subset of $Y$. We obtain
\be\label{eq:Zcapsupp} \del
\pi(Z_{P_\phi}^\delta) \cap \operatorname{supp} \phi = \emptyset
\ee
from (\ref{eq:<1/30}) and the choice $\delta = \frac{1}{10}
\|P_\phi\|_{C^0}$.

Then we choose the perturbation $\e$ of $P$ as given Lemma \ref{epsilon} so that
the inclusion
$$
Z_\e: = \{(\theta,y) \mid \epsilon (\theta,y) = 0 \} \subset \{(\theta,y) \mid
|P_\phi(y)| < \delta \} = Z_{P_\phi}^\delta
$$
holds. The function $h^{od}_\phi(\theta,y) - \frac{\del \widetilde
g}{\del t}(\theta,y) (= \e (\theta, y))$ is non-vanishing outside
$Z_\e$, in particular outside $Z_{P_\phi}^\delta$.

Now we closely describe the difference set
$$
\iota_{(h,\iota_Y)}^\phi \ominus dg.
$$
We would like to note that we have the base point formula
$$
\pi_2 \circ \pi(\iota_{(h,\iota_Y)}^\phi(\theta,y)) = y
$$
and
$$
\pi_1 \circ \pi(\iota_{(h,\iota_Y)}^\phi(\theta,y)) = \theta
$$
since $\phi$ is fiberwise-like, i.e., $\phi(\theta,y) =
(\theta,Y_\theta(y))$.

Therefore the projection in the direction $T^*S^1$ of
$\iota_{(h,\iota_Y)}^\phi $ consists of the elements of the form
\bea\label{eq:gthetay}
\alpha^\phi_g(\theta,y) & = & \widetilde \pi_1\circ
(\iota_{(h,\iota_Y)}^\phi)(\theta,y) - \left(\frac{\del g}{\del \theta}(\theta,y)
\right) d\theta \label{eq:alphaNphi} \\
& = &\alpha^\phi(\theta,y) - \left(\frac{\del g}{\del
\theta}(\theta,y)
\right) d\theta \nonumber \\
& = & \left(h^{od}_\phi(\theta,y)- \frac{\del g}{\del \theta}(\theta,y)
\right) d\theta \nonumber \\
& = & \e(\theta,y)\, d\theta
\eea

Here we use Lemma
\ref{hiphi} for the second equality, and the identity
\beastar
\e(\theta,y) = h^{od}_\phi(\theta,y)
- \frac{\del g}{\del \theta}(\theta, y).
\eeastar

On the other hand, $\overline U \setminus \pi(Z_{P_\phi}^\delta)$,
an open subset of $\overline U$ is disjoint from
$\partial U$ by Lemma \ref{Pphi}. Therefore Lemma
\ref{Morsek} enables us to find a smooth function $k$ on $\overline U$
that has no critical point on $\pi(Z_{P_\phi}^\delta) \supset \partial \overline U$ and
\be\label{dkn>0}
- dk(\vec n) > 0 \quad \mbox{on $\partial U$}.
\ee
By multiplying a large \emph{positive} constant,
if necessary, we may assume $dk$ to be arbitrarily large on
$\pi(Z_{P_\phi}^\delta)$ so that
\be\label{eq:dky} |dk(y)| > |\widetilde
\pi_2\circ (\iota_{(h,\iota_Y)}^\phi)(\theta,y) - d_y g(\theta, y)|
\ee
for $(\theta,y) \in Z_{P_\phi}^\delta \cap (S^1(2)
\times \overline U)$, and
\bea\label{eq:betakphi(n)>0}
\beta_k^\phi(\vec n) & = &\widetilde \pi_2(
(\iota_{(h,\iota_Y)}^\phi)(\theta,y))(\vec n)- \left(dk(y)(\vec n)
 + d_y g(\theta, y)(\vec n)\right)\nonumber \\
& = &- dk(y)(\vec n) + \left(\widetilde \pi_2\circ
(\iota_{(h,\iota_Y)}^\phi)(\theta,y)(\vec n)
 - d_y g(\theta, y)(\vec n)\right) > 0
\eea
on $\del(S^1(2)\times \overline U) = S^1(2)\times \partial \overline U$.

In particular (\ref{eq:dky}) implies \be\label{eq:pi2iotaphifk}
\beta_{k}^\phi(\theta,y) : = \widetilde \pi_2\circ
(\iota_{(h,\iota_Y)}^\phi)(\theta,y)- \left(dk(y) + d_y g(\theta, y)
\right) \neq 0 \ee for $(\theta,y) \in Z_{P_\phi}^\delta$. And
(\ref{eq:betakphi(n)>0}) implies that the projection of the image of
$\iota_{(h,\iota_Y)}^\phi \ominus d(g + k)$ to $T^*Y$ does not
intersect $\nu_-^*(\partial U)$.

On the other hand, on $\overline U \setminus Z_{P_\phi}^\delta$,
$\alpha_\e^\phi = \widetilde \pi_1 \circ \iota_{(h,\iota_Y)}^\phi(\theta,y)$
does not vanish by definition of $Z_{P_\phi}^\delta$.
Combining this with (\ref{eq:betakphi(n)>0})
and (\ref{eq:pi2iotaphifk}), we conclude that the difference set
$\iota_{(h,\iota_Y)}^\phi \ominus d(g +k)$ does not intersect
$$
\nu^*(S^1(2) \times \overline U)
$$
on $\overline U \cap Z_{P_\phi}^\delta$. Therefore we have proved that
$$
\iota_{(h,\iota_Y)}^\phi \ominus df
$$
does not intersect $\nu^*(S^1(2) \times \overline U)$, if we choose
$f = g + k$.

Finally to extend the function $f: \overline U \to \R$ to
$Y$, we first extend the function $f$ to a slightly bigger
open subset $S^1(2) \times U' \supset S^1(2) \times \overline U$
and then apply the partitions of unity subordinate to the
covering $\{U', Y \setminus \overline U\}$.
This finishes the proof.
\end{proof}

\medskip

\begin{proof}[Wrap-up of the proof of Theorem \ref{local-uniqueness}]
We set
\beastar
(X,\omega) & =& (M\times M, (-\omega_M)\oplus \omega_M), \quad Y = M\\
h_i(t,x) & = & H_i(t,\phi_{H_i}^t(x)), \quad \psi_i(t,x) =
(t,x, \phi_{H_i}^t(x)). \eeastar Then the hypothesis (1) of
Theorem \ref{local-uniqueness} implies $\psi_i$ uniformly
converges to the embedding $\iota_Y$ on $U$. On the other hand,
the hypotheses (1) and (2) therein, combined with the hypothesis that
$H_i$ is Cauchy in $C^0$-topology, imply
\be\label{eq:TanHiHi} \|\operatorname{Tan}(H_i)
- H_i\|_U \to 0 ; \quad \operatorname{Tan}(H)(t,x) :=
H(t,\phi_H^t(x)).
\ee
as $i \to \infty$, and so
$h_i=\operatorname{Tan}(H_i)$ is Cauchy on $[0,1] \times U$ as $H_i$
is assumed so. Therefore the pair $(h_i,\psi_i): [0,1] \times
\Delta \to [0,1] \times M\times M$ defines an exact Lagrangian
isotopy that satisfies all the hypotheses in Theorem
\ref{lag-localuniqueness}, except the normalization condition
(\ref{eq:average0}). However this can be easily adjusted by adding
a function $c_i$ on $[0,1]$ to $h_i$ and considering $H_i - c_i$
instead. This does not affect $\psi_i$ at all and so $(h_i-c_i,
\psi_i)$ is still an exact Lagrangian isotopy. From this, we
conclude $h-c_\infty \equiv 0$ and so $h\equiv c_\infty$. Hence
$h(t,x,x) \equiv c_\infty(t)$ for all $x \in U$. This in turn
proves that we have
\be\label{eq:TanHi0}
\|\operatorname{Tan}(H_i)- c_\infty \| \to 0.
\ee
Combining
(\ref{eq:TanHiHi}), (\ref{eq:TanHi0}), we have proved $H\equiv c_\infty$
on $[0,1] \times U$, and finished the proof of Theorem
\ref{local-uniqueness}.
\end{proof}

\section{Proof of Lemma \ref{notclosed}}
\label{sec:notclosed}

In this section, we prove Lemma \ref{notclosed}. Because the proofs
are not very different for the $C^0$ and $L^{(1,\infty)}$ cases, we
prove this closedness for more general class of
$L^{(1,\infty)}$-functions $h$ for the future study of
$L^{(1,\infty)}$-Hamiltonian flows.

Since $h \in L^{(1,\infty)} \subset L^1$, it follows that $\alpha = h^{od} \,d\theta$
defines a well-defined $n$-current on $S^1(2) \times Y$ (see
\cite{federer}). For the simplicity, we assume that $Y$ is
orientable and fix a volume form $\Omega$ on $Y$. For the
non-orientable case, exactly the same proof works by replacing the
real-valued current by one with coefficients in a flat line bundle
(or in the orientation sheaf) \cite{derham}.

It will be enough to find a smooth $(n-1)$-form $\eta$ on $S^1(2)
\times Y$ such that we have
\be\label{eq:alphawedgedeta}
\int_{S^1(2) \times Y} \alpha \wedge d\eta \neq 0
\ee
with $\operatorname{supp}\eta \subset U$. For the
$(n-1)$ form $\eta$, we can write
$$
d\eta = d\theta \wedge \frac{\del \eta}{\del \theta} + d_Y \eta
$$
where $d_Y$ is the exterior derivative in the direction of $Y$ and
$\frac{\del \eta}{\del \theta}: = \frac{\del}{\del\theta}\rfloor
d\eta$. We can write \be\label{eq:dYeta} d_Y\eta = f(\theta,y)
\Omega \ee for some smooth function $f$ on $S^1(2) \times U$ which
must satisfy \be\label{eq:intf=0} \int_U f(\theta,y) \, \Omega = 0
\ee for all $\theta \in S^1$. By applying a simple argument from the
Hodge theory, we can choose a family of $(n-1)$-forms $\eta_\theta$
on $U$ that satisfy (\ref{eq:dYeta}) for each $\theta\in S^1$ and
depend smoothly on $\theta$.

We next note that since $\alpha$ has the form $\alpha = h^{od}
d\theta$ we have
$$
\alpha \wedge d\eta = \alpha \wedge d_Y \eta
$$
and so derive
$$
\int_{S^1(2) \times U} \alpha \wedge d\eta = \int_{S^1(2) \times U}
\alpha \wedge (f \Omega) = \int_{S^1(2) \times U} f h^{od}\, d\theta
\wedge \Omega.
$$
Therefore it is enough to choose $\eta$ so that this integral does
not vanish. We will use $\eta$ of the form given by
$$
\eta = \pi_2^*(\eta_\theta)
$$
where $\eta_\theta$ are those chosen above. Now it remains to choose
$f$ above so that
\be\label{eq:intfhnot0} \int_{S^1(2) \times U} f
h^{od}\, d\theta \wedge \Omega \neq 0.
\ee
Since we have
$$
0 < \|h\| = \int_0^1 \text{osc}(h_t)\, dt
$$
by the hypothesis, there exists a measurable subset $A \subset [0,1]
\subset S^1$ with $0 < m(A)$ such that $\text{osc}(h_t)$ is defined
on $A \times U$ and $\text{osc}(h_t) > 0$ for $t\in A$. Furthermore
it is not difficult to check (see Theorem 2.6 in the original version of
\cite{oh:hameo2} for its proof) that $h_t$ is continuous on $U$ at such
$t \in A$. We consider the sets
\beastar
U_{A,h}^+ &:= &\{ (t,y) \in A \times U \mid h(t,y) > 0 \} \\
U_{A,h}^- &:= &\{ (t,y) \in A \times U \mid h(t,y) < 0 \}. \eeastar
Note that because of the normalized condition (\ref{eq:average0})
and $\text{osc}(h_t) > 0$ for $t \in A$, the sections
$$
U_{A,h}^\pm(t): = U_{A,h}^\pm \cap (\{t\} \times U)
$$
are non-empty open subsets of $U$ and so have positive measure. Then
it follows from the Fubini theorem that both $U_{A,h}^\pm$ have
positive measure, say $m_\pm > 0$ respectively. By definition,
$U_{A,h}^+ \cap U_{A,h}^- =\emptyset$. Furthermore, we also have
\bea
e^+_{A,h} & := & \int_{U^+_{A,h}} h \, d\theta \wedge \Omega > 0 \\
e^-_{A,h} & := & \int_{U^-_{A,h}} (- h) \, d\theta \wedge \Omega > 0
\eea and \bea
e^+_{A,h}(t) & := & \int_{U^+_{A,h}(t)} h_t \, \Omega > 0 \\
e^-_{A,h}(t) & := & \int_{U^-_{A,h}(t)} (- h_t) \,\Omega > 0 \eea
for $t \in A$.

By the standard outer-measure property of the measure induced by the
volume form $d\theta \wedge \Omega$, for any given $\e > 0$, we can
find open subsets $V^\pm_\e \subset [0,1] \times U$ such that
$$
V^\pm_\e \supset U_{A,h}^\pm
$$
and
$$
m(V^\pm_\e) \leq m(U_{A,h}^\pm) + \e, \quad m(V^+_\e \cap V^-_\e)
\leq \e.
$$
In particular we have \be\label{eq:mVUAe}
\begin{split}
m(V^+_\e\setminus \overline{V^-_\e}) &\geq m(U_{A,h}^+) - 2\e = m^+ - 2\e\\
m(V^-_\e \setminus \overline{V^+_\E}) &\geq m(U_{A,h}^-) - 2\e = m^-
- 2\e
\end{split}
\ee By choosing $\e$ so that $\e = \frac{1}{4}\min\{m^+,m^-\}$, it
follows that both $V^+_\e\setminus \overline{V^-_\e}$ and
$V^+_\e\setminus \overline{V^+_\E}$ are two disjoint open subsets of
positive measure.

Now we choose a partitions of unity $\{\chi_\e^+,
\chi_\e^-,\chi_\e^0\}$ subordinate to the open covering
$$
\{V^+_\e\setminus \overline{V^-_\e},\,  V^-_\e\setminus
\overline{V^+_\e}, \, M \setminus \overline{V^+_\e \cup V^-_\e}\}.
$$
We note that as $\e \to 0$, we have
$$
\chi_\e^\pm \to \chi_{U^\pm_{A,h}} \quad \mbox{in }\, L^1
$$
respectively, where $\chi_B$ denotes the characteristic function of
the set $B$ in general.

We consider the function $f_\e$ of the form
$$
f_\e(\theta,y) = a^+_\e(\theta,y) \chi_\e^+(\theta,y) -
a^-_\e(\theta,y) \chi_\e^-(\theta,y)
$$
for a suitable choice of positive functions $a^\pm_\e>0$ so that
\be\label{eq:intfe=0} \int_{\{\theta\} \times U} f_\e(\theta,\cdot)
\Omega = 0 \ee for all $\theta \in S^1$. We can choose
\be\label{eq:a-e} 0< a^\pm_\e(\theta,y) \leq C \ee for some $C > 0$
independent of $\e$ : We have only to define
\bea\label{eq:apme}
\begin{split}
a^+_\e(\theta) & = & \min\left\{1, \frac{\int_U
\chi_\e^-(\theta,y)\Omega}
{\int_U \chi_\e^+(\theta,y)\Omega} \right \}\\
a^-_\e(\theta) & = & a^+_\e(\theta) \left(\frac{\int_U
\chi_\e^+(\theta,y)\Omega} {\int_U
\chi_\e^-(\theta,y)\Omega}\right).
\end{split}
\eea
Then we note that as $\e \to 0$, the quotient
$\frac{\int_U
\chi_\e^-(\theta,y)\Omega}
{\int_U \chi_\e^+(\theta,y)\Omega}$
converges to $\frac{m(U^-_{A,h})\setminus m(U^+_{A,h})}
{m(U^+_{A,h})\setminus m(U^-_{A,h})}$ and hence
\bea
0 & < & a_\e^+(\theta) \leq \frac{m(U^-_{A,h})\setminus m(U^+_{A,h})}
{m(U^+_{A,h})\setminus m(U^-_{A,h})} + 1 \\
0 & < & a_\e^-(\theta) \leq \frac{m(U^-_{A,h})\setminus m(U^+_{A,h})}
{m(U^+_{A,h})\setminus m(U^-_{A,h})} + 1.
\eea
We just choose $C$ to be
$$
C = \max\left\{\frac{m(U^-_{A,h})\setminus m(U^+_{A,h})}
{m(U^+_{A,h})\setminus m(U^-_{A,h})} + 1,\frac{m(U^-_{A,h})\setminus m(U^+_{A,h})}
{m(U^+_{A,h})\setminus m(U^-_{A,h})} + 1\right\}.
$$
This in turn implies \beastar \int_{S^1\times U}\chi_A(\theta)
f_\e h^{od}d\theta \wedge \Omega
& = & \int_{[0,1] \times U} \chi_A(\theta) f_\e h\, d\theta \wedge \Omega \\
& \to & \int_{[0,1] \times U} h(a_{(+,A)}\chi_{U^+_{A,h}} -
a_{(-,A)}\cdot\chi_{U^-_{A,h}}) d\theta \wedge \Omega \eeastar as
$\e \to 0$, where we denote by $a_{(\pm,A)}$ the $L^1$-limits of
$a^\pm_\e$. The last integral is strictly positive because both
$a_{(\pm,A)}$ cannot vanish simultaneously by the choice of
$a^\pm_e$ made in (\ref{eq:apme}).

Therefore if $\e > 0$ is sufficiently small,
(\ref{eq:alphawedgedeta}) holds for $\eta$ given by $\eta =
\pi_Y^*(\eta_\theta)$ where $\eta_\theta$ is the form satisfying
(\ref{eq:dYeta}) for $\varphi_\e(\theta) f_\e $ with $\varphi_\e =
\varphi_\e(\theta)$ being a suitable smooth $L^1$-approximation of
the characteristic function $\chi_A$ on $S^1$. This proves that the
current $\alpha$ is not closed and so finishes the proof of Lemma
\ref{notclosed}.

\section{Proof of Lemma \ref{alphaphi} : Problem of mass transport}
\label{sec:transport}

Let $\alpha = h\, d\theta$ as in section \ref{sec:analysis}.
Again we will prove this for $h$ lying in $L^{(1,\infty)}$. (In
this section, for the simplicity of notations, we will just denote
$h$ for $h^{od}$.)

We re-state  here in its contrapositive form,
which is more close to the statement in Lemma 2.5 \cite{viterbo2}.

\begin{lem}\label{viterbolemma} Let $\alpha$ be as above.
Suppose
$$
\int_{\pi_Y^{-1}(y)} \phi_*\alpha = 0
$$
for all $y \in U$ and for any fiberwise-like diffeomorphism
$\phi:S^1(2) \times Y \to S^1(2) \times Y$ with
$\operatorname{supp}\phi \subset U$. Then $\alpha$ is closed
as an $n$-current on $S^1(2) \times U$.
\end{lem}

We recall that for any fiberwise-like diffeomorphism the map
$Y_\theta: Y \to Y$ is a diffeomorphism for each $\theta \in S^1$,
where $Y_\theta(y) = Y(\theta,y)$ in the representation of $\phi$,
$$
\phi(\theta,y) = (\theta,Y(\theta,y)).
$$
The rest of the section will be occupied with the proof of this
lemma.

Let $\eta'$ be an $(n-1)$ form on $U$ and consider the
pull-back form $\pi^*\eta'$ on $S^1(2) \times U$. Then the Fubini
theorem implies
$$
\int_{S^1(2) \times U} \phi_*\alpha \wedge \pi^*d\eta' = \int_U
\left(\int_{S^1} \phi_*\alpha\right) \wedge d_Y\eta' = 0.
$$
Therefore we have \beastar 0 & = &\int_{S^1(2) \times U} \alpha
\wedge \phi^*\pi^*(d\eta')
=\int_{S^1(2) \times U} \alpha \wedge (\pi\circ \phi)^*d\eta' \\
& = & \int_{S^1(2) \times U} \alpha \wedge Y_\theta^*(d\eta') =
\int_{S^1(2) \times U} \alpha \wedge d(Y_\theta^*\eta') \eeastar Now
the proof of Lemma \ref{viterbolemma} will be finished by the
following proposition

\begin{prop} Let $\alpha = h\, d\theta$ be the above
section with coefficients in $L^{(1,\infty)}$ considered as an $n$-current
considered above on $S^1(2) \times U$. Suppose that
\be\label{eq:alphaeta'} \int_{S^1(2) \times U} \alpha \wedge
d(Y_\theta^*\eta') = 0 \ee for all smooth $(n-1)$ forms $\eta'$ on
$U$ and for all fiberwise-like diffeomorphism $\phi: S^1(2) \times
Y \to S^1(2) \times Y$ with $\phi(\theta,y) = (\theta,
Y(\theta,y))$, where $Y_\theta: Y\to Y$ is the diffeomorphism $y
\mapsto Y(\theta,y)$ such that $\operatorname{supp}Y_\theta \subset
U$ for all $\theta \in S^1(2)$.
Then $\alpha$ is closed on $S^1(2) \times U$.
\end{prop}
\begin{proof}
We have to show
\be\label{eq:alphaclosed} \int_{S^1(2) \times Y}
\alpha \wedge d\eta = 0
\ee for
all smooth $(n-1)$-form $\eta$ supported in $S^1(2) \times U$.

Here we note that closedness (\ref{eq:alphaclosed})
of a current is a local property and so
it is enough to check the closedness in a coordinate neighborhood
$(a,b) \times B \subset S^1(2) \times U$ of a given point $(e,y) \in
S^1(2) \times U$. By restricting to an even smaller neighborhood
$(c,d) \times C \subset \overline [c,d] \times C \subset (a,b)
\times B$ of point $(e,y)$, we may assume that $d\eta$ is supported
on $(c,d) \times C \cup (Y \setminus \overline B)$.
We will assume this for the rest of the proof.

We decompose
$$
\eta = \eta_Y + d\theta \wedge \zeta_Y
$$
where $\eta_Y$ and $\zeta_Y$ are forms satisfying
$$
\frac{\del}{\del\theta} \rfloor \eta_Y = \frac{\del}{\del\theta}
\rfloor \zeta_Y =0.
$$
Since $\alpha = h\, d\theta$, it follows
$$
\alpha \wedge d\eta = \alpha \wedge d_Y(\eta_Y).
$$
Now we would like to compare $d_Y(Y_\theta^*(\eta'))$ and
$d_Y(\eta_Y)$. Both are $n$-forms supported in $S^1(2) \times U$ with coefficients depending
on $(\theta,y)$. More precisely, we have
$$
d_Y(\eta_Y) = f(\theta,y)\Omega, \quad \int_Y  f(\theta,y)\Omega = 0
$$
and
$$
d_Y(Y_\theta^*(\eta')) = Y_\theta^*(d_Y(\eta')) =
Y_\theta^*(g\Omega), \quad \int_Y  g(y)\Omega = 0.
$$
For given $f$ as above, we would like to represent
$f(\theta,y)\Omega$ as $Y_\theta^*(g\Omega)$, i.e., we would like to
solve the equation \be\label{eq:fgY} f(\theta,y)\Omega =
Y_\theta^*(g\Omega)\quad  \Big( = g(Y(\theta,y))|\det
(dY_\theta)|\Omega\Big) \ee for a suitable choice of $\phi$ (and so
$Y_\theta$) and $g:Y \to \R$.

We will apply Moser's trick to solve this equation. Since we can
make the coefficients at a fixed point $\theta$ the same on both sides, it
is enough to consider its derivative.  We write (\ref{eq:fgY}) as
$$
(Y_\theta)_*\left(f(\theta,\cdot)\Omega\right) = g(y)\Omega.
$$
Differentiating this with respect to $\theta$, we get
$$
(Y_\theta)_*\left(\frac{\del f}{\del \theta}\Omega -d_Y(X_\theta
\rfloor f \Omega\right) =0
$$
where $X_\theta$ is the vector field generating the family
$Y_\theta: Y \to Y$ of diffeomorphisms. Since $Y_\theta$ is a
diffeomorphism supported in $U$, we obtain
\be\label{eq:dfdtheta} \frac{\del f}{\del
\theta}\Omega - d_Y(X_\theta \rfloor f \Omega) = 0
\ee
with $\operatorname{supp} X_\theta \subset U$. Let
$\beta_\theta$ be an $(n-1)$ form on $Y$ such that
\be\label{eq:dbeta} \frac{\del f}{\del \theta}\Omega = d_Y
\beta_\theta \ee which exists because
$$
\int_Y \frac{\del f}{\del \theta}\Omega = \frac{\del}{\del \theta}
\int_Y  f\, \Omega = 0.
$$
We may assume that $\beta_\theta$ depends smoothly on $\theta$.
Substituting (\ref{eq:dbeta}) into (\ref{eq:dfdtheta}), the problem
is reduced to solving the algebraic equation \be\label{eq:beta}
\beta_\theta = f X_\theta \rfloor \Omega \ee in terms of $X_\theta$.
This equation can be solved pointwise whenever $f$ does not vanish.

By adjusting
$f\Omega$ outside $[a,b] \times B$ so that $\int_Y f \Omega = 0$ is
satisfied, we may safely assume that
\be\label{eq:fneq0}
f \neq 0 \quad \mbox{ on $(c,d) \times C$}.
\ee
Then we can solve (\ref{eq:beta}) for $X_\theta$ on $(c,d) \times B$
which will have support on $(c,d) \times C \cup (U \setminus
\overline B)$.
Now by making the interval $(c,d)$ smaller if
necessary, we can integrate
$$
\dot y = X_\theta(y), \quad y(e) = y
$$
to obtain a flow $Y_\theta$ that solves (\ref{eq:fgY}) on $(c,d)
\times C$, if we choose the function $g$ above so that its values
are given by
$$
g(y) = f(e,y)|\det(dY_e(y))|.
$$
Therefore the hypothesis (\ref{eq:alphaeta'}) implies that $\alpha$
is closed on $(c,d) \times C$. Since this holds at any given point
$(e,y) \in S^1(2) \times U$, $\alpha$ is a closed current on $S^1(2)
\times U$. This
finishes the proof of Lemma \ref{viterbolemma}. Now the statement in
Lemma \ref{alphaphi} concerning the $C^0$-closeness to the identity
is immediate from the proof, because $\phi$ can be obtained by
integrating a vector field $X_\theta$ on a short `time' interval
starting from $\theta = e$ in this proof.

The above also proves that the diffeomorphism $\phi$ constructed
above has the form $\phi(\theta,y) = (\theta, Y_\theta(y))$.
Hence the proof of Lemma \ref{alphaphi}.
\end{proof}

\section{The case of general boundary and etc.}
\label{sec:general}

In this section, we explain how we derive the same result for
a general open subset whose boundary is not necessarily smooth
and how to treat general $C^0$-Hamiltonian path which is
not necessarily boundary-flat. We then explain how we can
generalize to open subsets $U \subset M$ where $M$
is not necessarily closed. Finally we explain how the proof
of the uniqueness of compactly supported Hamiltonians, Theorem \ref{uniqueness},
either on $M$ with contact-type boundary or on $M$ which
is bounded at infinity.

For the first point, we just choose an increasing sequence $U_l \subset U$
of open subsets
such that
$$
U = \bigcup_{l} U_l
$$
and each $U_l$ has smooth boundary $\partial U_l$. Obviously if
$\phi_{H} = \lambda \equiv id$ on $U$, so is on $U_i$ and so we
obtain
$$
H(t,x) \equiv c_l(t)
$$
on $U_l$ for each $i$. Since $U_l \subset U_{l'}$ for $l' \geq l$,
we must have $c_l \equiv c_{l'}$ for all $l, \, l'$. Defining $c$ to
be the common function, we have proved $H(t,x) \equiv c(t)$ on $U$
and so on $\overline U$ by continuity of $H$.

For the second point, we imitate the idea of Viterbo from [Erratum, \cite{viterbo2}].
Let $H_i$ be a sequence of Hamiltonians, which are not
necessarily boundary-flat, such that $\phi_{H_i}$ and $H_i$
converges in $C^0$-topology.

Then for any reparameterization $H_i^{\zeta}$
defined by $H^\zeta(t,x) = \zeta'(t)H(\zeta(t),x)$ with each
fixed function $\zeta: [0,1] \to [0,1]$, both
$\phi_{H_i^\zeta}$ and $H_i^\zeta$ converge in $C^0$-topology.
Furthermore $\phi_{H_i^\zeta}$ converge to $id$ since
$\phi_{H_i^\zeta}$ is nothing but the time reparameterization
$\phi_{H_i}^{\zeta(t)}$ and $\zeta$ fixed. Therefore
we have
$$
H^\zeta \equiv c_\zeta \quad \mbox{on $U$}.
$$
We note that for any given $t_0 \in [0,1]$ we can choose some
$\zeta_{t_0}$ so that $\zeta_{t_0}(t_0) = t_0$ and $\zeta'_{t_0}(t_0) = 1$.
By considering such $\zeta$, we prove
that
$$
H(t_0,\cdot) = c_{\zeta_{t_0}}(t_0)
$$
and hence $H(t,x) \equiv c_{\zeta_t}(t)$. We just define $c:[0,1]
\to [0,1]$ to be
$$
c(t) : = c_{\zeta_t}(t).
$$
Since $H$ is continuous, this $c$ is a continuous function on
$[0,1]$. This takes care of the case of arbitrary
$C^0$-Hamiltonian paths and finishes the proof of Theorem
\ref{local-uniqueness}.

For the point of generalizing the locality to an arbitrary complete symplectic manifold
\emph{bounded at infinity} in the sense of Gromov \cite{gromov:pseudo} but
not necessarily closed, we have only to point out that
the proof of the relative version in \cite{kasturi-oh} of Lagrangian
intersection result between the Hamiltonian deformation of
the zero section $o_{S^1(2) \times M}$ with the conormal
$$
\nu^*(S^1(2) \times \overline U)
$$
goes through verbatim, once we have an a-priori $C^0$-estimate for
the relevant Floer trajectory moduli space, which makes
the Floer homology between $\nu^*(S^1(2) \times \overline U)$
and the zero section $o_{S^1(2) \times Y}$ well-defined
and invariant under the compactly supported Hamiltonian
isotopy of $o_{S^1(2) \times Y}$
with $Y = M$ by the standard continuation argument.
Here we also use the fact that the intersection
$$
\nu^*(S^1(2) \times \overline U) \cap o_{S^1(2) \times Y}
$$
is compact. We leave the details to interested readers.

An immediate consequence of this locality is the following uniqueness
theorem of compactly supported topological Hamiltonians introduced
in \cite{oh:hameo1}.

\begin{thm}\label{compactsupp}
Let $M$ be a symplectic manifold which is bounded at infinity.
Suppose $H_i:[0,1] \times M \to \R$ is a sequence of smooth
normalized Hamiltonian functions such that
\begin{enumerate}
\item there exists a compact set $K \subset \operatorname{Int}M$
such that
$$
\operatorname{supp} H_i \subset K
$$
\item $H_i$ converges to $H: Y \to \R$ in $C^0$-topology.
\end{enumerate}
Then if $\phi_{H_i} \to id$ uniformly on $[0,1] \times M$, we have
$H \equiv 0$.
\end{thm}
\begin{proof} We choose $U$ so that $\del U$ is smooth, $\overline U$
compact and $K \subset U$. Then we apply the locality theorem.
\end{proof}

Finally the uniqueness result for compact $(M,\omega)$ with
contact-type boundary $\del M$ is reduced to the above case of
complete manifold bounded at infinity by attaching the end $\del M
\times \R_+$ along its boundary. This is standard and so left to
readers for the details. It is not clear to the author whether
or not the uniqueness still holds for general $(M,\omega)$ with
smooth boundary, not necessarily of contact-type, for which
there does not seem to exist an obvious way of symplectically
embedding $(M,\omega)$ into a symplectic manifold bounded at infinity.

\section{Appendix}
\label{sec:appendix}

\subsection{Exact Lagrangian property of the isotopy (\ref{eq:f})}
\label{subsec:exactproperty}

In this appendix, we prove the exactness of the Lagrangian embedding
$\iota_{(h,\psi)}$ given in (\ref{eq:iotahpsi}). This is a
consequence of Theorem 6.1.B \cite{polterov} whose proof is left as
an exercise in \cite{polterov} with the indication that the
exactness of $i_s^*\alpha$ will follow from some homological
argument. Here we give a computational proof somewhat different from
the proof suggested in \cite{polterov} and also explicitly write
down the anti-derivative of $i_s^*\alpha$. The outcome is quite
interesting and suggestive, and gives rise to a slight
generalization of Theorem 6.1.B \cite{polterov} in that it applies
to any Hamiltonian paths, not just to Hamiltonian loops, to which
homological argument may not apply.

Consider an isotopy of Hamiltonian paths i.e., a two parameter
family of Hamiltonian diffeomorphisms
$$
(s,\theta) \in [0,1] \times [a,b] \mapsto \psi_s^\theta \in Ham(X,\omega).
$$
We denote by $H =H(s,\theta,x)$ the \emph{normalized} Hamiltonian generating the $\theta$-isotopy
$\theta \mapsto \psi_s^\theta$ and $K = K(s,\theta,x)$ the one generating
the $s$-isotopy $s \mapsto \psi_s^\theta$. In other words,
$$
\frac{\del \psi^\theta_s}{\del \theta}\circ (\psi^\theta_s)^{-1} = X_H,\quad
\frac{\del \psi^\theta_s}{\del s}\circ (\psi^\theta_s)^{-1} = X_K.
$$
Then the following identity has been proved in \cite{banyaga78},
\cite{polterov}, \cite{oh:jkms} \be\label{eq:flat} \frac{\del
K}{\del \theta} = \frac{\del H}{\del s} - \{K,H\} \ee for the
normalized Hamiltonian $H$.

With this preparation, we prove the following theorem.

\begin{thm}[Compare with Theorem 6.1.B \cite{polterov}]\label{exactisotopy}
Let $Y \subset X$ be a Lagrangian embedding and consider a two
parameter family $\psi_s^t$ of Hamiltonian diffeomorphisms as above.
Then the isotopy $f: [0,1] \times [a,b] \times Y \to T^*[a,b] \times
X=: M $ defined by
\be\label{eq:f}
f(s, \theta, x) = (\theta, - H(s,
\theta,\psi^\theta_s(x)),\psi^\theta_s(x))
\ee
is an exact Lagrangian isotopy in the sense of Definition \ref{exact-lag}.
Furthermore, we have
\be\label{eq:isalpha}
i_s^*\alpha = d(-K\circ f)_s.
\ee
\end{thm}

We like to mention that here $s$ replaces the role of $t$
and $[a,b] \times Y$ that of $Y$ in Definition \ref{exact-lag}.

\begin{proof}
First we compute
$$
f^*(d\theta \wedge db + \omega_X).
$$
Since $o_{[a,b]} \times i_Y \subset T^*[a,b] \times X$
is a Lagrangian embedding and $\psi_s^\theta$ is symplectic, it is easy to check that
$$
f^*(d\theta \wedge db + \omega_X) = \alpha \wedge ds
$$
for some one-form $\alpha$ on $[0,1] \times ([a,b] \times Y)$.
We have
$$
\alpha = \frac{\del}{\del s} \rfloor f^*(d\theta \wedge db + \omega_X)
$$
which we now evaluate. First we have
$$
\alpha\left(\xi \right) = (d\theta \wedge db + \omega_X)
\left(\frac{\del f}{\del s}, Tf (\xi)\right)
$$
for general tangent vector $\xi \in T([a,b] \times Y)$. Straight computations
give rise to the following formulae :
\beastar
Tf \left(\frac{\del}{\del s}\right) & = &
 \left (\frac{\del H}{\del s} - \{K,H\}\right) \frac{\del}{\del b}
\bigoplus X_K\circ \psi^\theta_s\\
& = & \frac{\del K}{\del \theta} \frac{\del}{\del b} \bigoplus X_K\circ \psi^\theta_s\\
Tf \left(\frac{\del}{\del \theta}\right) & = & \frac{\del}{\del \theta}
\bigoplus - \frac{\del H}{\del \theta}
\frac{\del}{\del b} \bigoplus X_H\circ \psi^\theta_s\\
Tf(\eta)& = & -d_X H (T\psi^\theta_s(\eta))\frac{\del}{\del b}
\bigoplus T\psi^\theta_s(\eta)
\eeastar
where $\eta$ is a vector field whose projection to $[0,1] \times [a,b]$
is zero. For the first equality, we used the identity
$dH(X_K) = \{H,K\}$ and for the second equality, we have used (\ref{eq:flat}).
From these formulae, we derive
\beastar
i_s^* \alpha & = &  i_s^*\left(\frac{\del}{\del s} \rfloor(d\theta \wedge db + \omega_X)\right)
= -\left(\left(\frac{\del K}{\del \theta}\circ f_s \right)\, d\theta \circ Tf_s +
d_X K \circ Tf_s \right)\\
& = & -f_s^*\left(\frac{\del K}{\del \theta}\, d\theta + d_X K\right)
= -f_s^*dK = d(-K \circ f)_s
\eeastar
Therefore the form $i_s^*\alpha$ is exact
and hence $f$ defines a Lagrangian isotopy in the sense of Definition
\ref{exact-lag}. Furthermore we have also proved (\ref{eq:isalpha}).

This finishes the proof.
\end{proof}

To apply this theorem to the embedding $\iota_{(h,\psi)}$,
we replace $[a,b]$ by $[0,2]$
replace $[a,b]$ by $[0,2]$ and apply the theorem to $\psi^\theta_s
= \phi_H^{s\chi(\theta)}$ where $\chi:[0,2] \to [0,1]$ is a function
defined by
\be\label{eq:chi}
\chi(s) = \begin{cases} \zeta(s) & \quad \mbox{for } \, s \in [0,1] \\
\zeta(2-s) & \quad \mbox{for } \, s \in [1,2].
\end{cases}
\ee
Here $\zeta:[0,1] \to [0,1]$ is a reparameterization function as defined
in (\ref{eq:zeta}).

Recall the Hamiltonian generating
in $\theta$-direction, which we still denote by $H$ with some abuse of
notations, is given by
$$
H(s,\theta,x) = s\chi'(\theta)H(s\chi(\theta),x).
$$
Now we have only to check
that the Hamiltonian $K: [0,1] \times [0,2] \times Y \to \R$
generating in $s$-direction above is 2-periodic, i.e., satisfies
\be\label{eq:0=2}
K(s, 0, x) = K(s,2,x)
\ee
for the present case. But a straightforward calculation
differentiating $\psi_s^\theta$ in $s$ gives rise to
$$
K(s,\theta,x) = \chi(\theta) H(s\chi(\theta),x)
$$
which obviously satisfies (\ref{eq:0=2}) by the definition
(\ref{eq:chi}) of $\chi$.

\subsection{Proof of Lemma \ref{epsilon}}

In this appendix, we prove Lemma \ref{epsilon}.
We restate the lemma here.

\begin{lem} There exists a continuous function
$\e: S^1(2) \times Y \to \R$ such that
\begin{enumerate}
\item For all $y \in Y$, $\int_{S^1(2)} (h_\phi^{od}(t,y) - \e(t,y)) \, dt = 0$,
\item $h_\phi^{od}-\e$ is smooth.
\item We have the inequality
$$
\|P_\phi - \e \|_{C^0} < \frac{1}{30}\|P_\phi\|_{C^0}.
$$
\end{enumerate}
\end{lem}
\begin{proof}
Consider the solution $\widetilde f$ for the equation
(\ref{eq:dfdt}), which can be written as
$$
\widetilde f(\theta,y) = \int_0^\theta(h_\phi^{od}(t,y) - P_\phi(y))
\, dt.
$$
From this expression, it follows that $\widetilde f$ is
$C^1$ in $\theta$ and $C^0$ in $y$. Therefore we can choose a smooth
approximation $f_{sm}: S^1(2) \times Y \to \R$ so that
both norms
\be\label{eq:norms}
\left|\frac{\del \widetilde f}{\del \theta} - \frac{\del f_{sm}}{\del \theta}
\right|_{C^0}, \, |\widetilde f - f_{sm}|_{C^0}
\ee
as small as we want. Denote $\xi =  \widetilde f -f_{sm}$, and define
\be\label{eq:epsilondef}
\e(\theta,y): = P_\phi(y) + \frac{\del \xi}{\del\theta}.
\ee
We now show that $\e$ will do our purpose. First (1) follows since
we have
$$
h_\phi^{od} - \e = (h_\phi^{od}-P_\phi) - \frac{\del \xi}{\del\theta}
$$
and $h-\e$ satisfies the vanishing of the corresponding integral.

For (2), we recall
$$
h_\phi^{od}- P - \frac{\del \xi}{\del\theta}
= \frac{\del \widetilde f}{\del \theta} - \frac{\del \xi}{\del\theta}
= \frac{\del}{\del \theta}(\widetilde f - \xi) = \frac{\del f_{sm}}{\del \theta}
$$
which is smooth.

For (3), we note
$$
\e - P_\phi = \frac{\del \xi}{\del\theta}
= \frac{\del}{\del\theta}(\widetilde f - f_{sm})
$$
which can be made as small as we want by choosing $f_{sm}$
so that the norms in (\ref{eq:norms}) become small.
This finishes the proof.
\end{proof}

\end{document}